\documentclass[11pt]{article}
\def\baselinestretch{1.1}
\usepackage{amsmath}
\usepackage{amsthm}
\usepackage{amsfonts}
\usepackage{amssymb}
\usepackage{times}
\newcommand{\beq}[1]{ 
 \begin{equation}\label{#1}}
\def\ot{\otimes}
\def\g{{\rm g}}
\def\tu{\tilde{u}}
\newcommand{\fr}[2]{{\textstyle \frac{#1}{#2} }}
\def\klein{\scriptscriptstyle}
\newcommand{\al}{\alpha}

\newcommand{\fsl}{\mathfrak{sl}}
\def\sz{{\mathsf z}}
\def\Uq{{\mathcal U}_q(\fsl_2)}
\def\UqR{{\mathcal U}_q(\fsl(2,\!\BR))}
\def\Aq{{\mathcal A}_q}
\def\Wq{{\mathcal W}_q}
\def\Sq{{S\!L_q(2)}}
\def\Gq{{G\!L_q(2)}}
\def\SqR{S\!L\vphantom{L}_q(2,\!\BR)}
\def\GqR{G\!L\vphantom{L}_q(2,\!\BR)}
\def\Gqe{\widetilde{G\!L}_q(2)}

\def\GqeR{\widetilde{G\!L}\vphantom{L}_q(2,\!\BR)}
\def\GqeRp{\widetilde{G\!L}\vphantom{L}'_q(2,\!\BR)}
\def\GqeRpp{\widetilde{G\!L}\vphantom{L}''_q(2,\!\BR)}
\newcommand{\nn}{\nonumber}
\DeclareMathOperator*{\qdet}{qdet}
\newcommand{\0}{{\mathfrak 0}}
\newcommand{\1}{{\mathfrak 1}}
\newcommand{\2}{{\mathfrak 2}}
\newcommand{\3}{{\mathfrak 3}}
\newcommand{\4}{{\mathfrak 4}}
\def\p1{{\scriptscriptstyle +1}}
\newcommand{\FU}{{\mathfrak U}}
\newcommand{\rf}[1]{(\ref{#1})}
\newcommand{\CQ}{{\mathcal Q}}
\newcommand{\CA}{{\mathcal A}}
\newcommand{\CB}{{\mathcal B}}
\newcommand{\CT}{{\mathcal T}}
\newcommand{\CK}{{\mathcal K}}

\newcommand{\la}{\lambda}
\newcommand{\SRN}{{\rm N}}

\def\QED{\phantom{x}\hfill{\small$\square$}\par\vspace{5pt}}

\newcommand{\kl}{{\rm\klein L}}

\newcommand{\BR}{{\mathbb R}}

\newcommand{\BC}{{\mathbb C}}

\newcommand{\CW}{{\mathcal W}}
\newcommand{\CH}{{\mathcal H}}
\newcommand{\SR}{{\mathsf R}}
\newcommand{\SF}{{\mathsf F}}
\newcommand{\SH}{{\mathsf H}}
\newcommand{\SP}{{\mathsf P}}
\newcommand{\ST}{{\mathsf T}}
\newcommand{\SX}{{\mathsf X}}
\newcommand{\sr}{{\mathsf r}}
\newcommand{\sw}{{\mathsf w}}

\theoremstyle{plain}
\newtheorem{thm}{Theorem}
\newtheorem{propn}{Proposition}
\newtheorem{lem}{Lemma}
\newtheorem{defn}{Definition}
\theoremstyle{remark}
\newtheorem{rem}{Remark}
\begin{document}
\thispagestyle{empty}
%
\begin{center}
{\large\bf
Baxterization of $\Gq$ and its
application to \\ the Liouville model
and some other models on a lattice
} \\ [3mm]
{\sc Andrei G. Bytsko}    \\ [2mm]
{ \small
   Steklov Mathematics Institute,
 Fontanka 27, 191023, St.~Petersburg, Russia \\
 DESY Theory Group, Notkestrasse 85,
 D--22603, Hamburg, Germany
} \\ [2.5mm]
{ }

\end{center}
\vspace{1mm}
\begin{abstract}
We develop the Baxterization approach to (an extension of)
the quantum group~$\Gq$. We introduce two matrices which play
the role of spectral parameter dependent L--matrices and
observe that they are naturally related to two different
comultiplications. Using these comultiplication structures,
we find the related fundamental R--operators in terms of powers
of coproducts and also give their equivalent forms in terms of
quantum dilogarithms. The corresponding quantum local
Hamiltonians are given in terms of logarithms of positive
operators. An analogous construction is developed for the
\hbox{q--oscillator} and Weyl algebras using that
their algebraic and coalgebraic structures can be obtained
as reductions of those for the quantum group.
As an application, the lattice Liouville model, the q--DST
model, the Volterra model, a lattice regularization of the
free field, and the relativistic Toda model are considered.
\end{abstract}

{\small
\tableofcontents
}

\section{Introduction: motivation and outline of main results}

\renewcommand{\%}{\empty}

A quantum model is a system $(\CH,\CA,\SH)$, with a Hilbert space
$\CH$, an algebra of observables $\CA$, and a Hamiltonian~$\SH$. The model
is {\em integrable} if there exists a complete set of quantum integrals of
motion, i.e., a set of self--adjoint elements of $\CA$ which commute with
each other and with the Hamiltonian. For {\em homogeneous} one--dimensional
lattice models one has $\CH=\CK^{\ot \SRN}$, $\CA=\CB^{\ot \SRN}$, with
one copy of Hilbert space $\CK$ and algebra of local observables $\CB$
being associated to each of the $\SRN$ sites of a one--dimensional
lattice. $\CK$ is usually characterized as a representation of an
algebra $\FU$ of ``symmetries'', and $\CB$ is generated from the
operators which represent the elements of $\FU$ on~$\CK$.

A key step in constructing an integrable lattice model is to find
an \hbox{L--matrix} \hbox{$L(\la) \in {\mathrm{Mat}}(n)\otimes\CB$}
and an auxiliary \hbox{R--matrix}
$R(\la) \in {\mathrm{Mat}}(n)^{\otimes 2}$
such that the following matrix commutation relation
\beq{RLL0}
 R_{\1\2}(\la) \, L_{\1\3} (\la\mu) \, L_{\2\3} (\mu) =
 L_{\2\3} (\mu) \, L_{\1\3} (\la\mu) \, R_{\1\2}(\la) \,,
\end{equation}
where $\la,\mu\in \BC$, is equivalent to the defining relations
of~$\FU$. Here and below we use the standard notation: subscripts
indicate nontrivial components in tensor product, e.g.,
$R_{\1\2} \equiv R \otimes {\sf 1}$, etc. For further details
on the \hbox{R--matrix} approach to quantum integrability we refer
the reader to the review~\cite{F1}.

For a model on a {\em closed} one--dimensional lattice, i.e.,
with periodic boundary conditions, a set of quantum integrals
of motion is generated by the {\em auxiliary} transfer--matrix
\hbox{$T(\la) = \mathrm{tr}_a
   \bigl(L_{a,\SRN}(\la) \ldots L_{a,\1}(\la)\bigr)$}.
However, these integrals are in general non--local, i.e., they
are not representable as a sum of terms each containing nontrivial
contributions only from several nearest sites.
 The recipe \cite{FT2} for constructing {\em local} integrals
of motion for a model with a given L--matrix
  is to find first the corresponding
{\em fundamental} R--operator $\SR(\la) \in \CB^{\otimes 2}$,
which satisfies the following intertwining relation
(here and below we will use it in the braid form):
\beq{RLL}
 \SR_{\2\3}(\la) \, L_{\1\2} (\la\mu) \, L_{\1\3} (\mu) =
  L_{\1\2} (\mu) \, L_{\1\3} (\la\mu) \, \SR_{\2\3}(\la) \,.
\end{equation}
The corresponding transfer--matrix is constructed as
\hbox{$\ST(\la) = \mathrm{tr}_a  \bigl(\SR_{a\SRN}(\la)\SP_{a\SRN}
 \ldots \SR_{a\1}(\la)\SP_{a\1}\bigr)$},
where the subscript $a$ stands now for an auxiliary copy of~$\CB$,
and $\SP$ is the unitary operator permuting tensor factors
in~$\CB^{\otimes 2}$. The fundamental R--operator is usually
{\em regular}, that is, after appropriate normalization,
it satisfies the relation
\beq{reg}
  \SR_{}(1) = {\sf 1} \ot {\sf 1} \,.
\end{equation}
If the regularity condition holds, then first and higher order
logarithmic derivatives of $\ST(\la)$ at $\la{=}1$ are
local integrals of motion for the periodic homogeneous model in question.
In particular, the Hamiltonian is often chosen as the most local
integral which involves only nearest neighbour interaction:
\beq{Hloc}
 \SH = i
 \frac{\partial}{\partial \la}  \bigl[ \log \ST(\la) \bigr]_{\la=1}
 = \sum_{n=1}^{\SRN} H_{n,n \p1} = \sum_{n=1}^{\SRN}
 i \frac{\partial}{\partial \la}
 \bigl[ \SR_{n \p1,n}(\la) \bigr]_{\la=1} \,,
\end{equation}
where the summation assumes that $\SRN{+}1\equiv 1$.

Thus, finding the fundamental \hbox{R--operator} for a given
\hbox{L--matrix} is an important part of the R--matrix approach
to quantum integrable models. Furthermore, this problem is
closely related to the problem of constructing the corresponding
evolution operators and Q--operators.
However, there is no general
method for solving equation~\rf{RLL}. The particular
difficulty here is that it is not clear apriori on which
operator argument(s) the function $\SR(\la)$ depends.

Among the few known examples of constructing a fundamental
R--operator the most algebraically transparent are those
related to the case where the symmetry $\FU$ admits the structure
of a bialgebra. Such examples include the XXX spin chain
\cite{KRS} and closely related nonlinear Schr\"odinger model
\cite{FT2}, where $\FU={\cal U}(\fsl_2)$; and the
XXZ spin chain \cite{J1} and closely related sine--Gordon model
\cite{FT2,T1}, where $\FU=\Uq$.
A crucial observation for solving \rf{RLL} in these cases
is that the operator argument of $\SR(\la)$  is $\Delta(C_q)$,
where $C_q$ is the Casimir element of $\FU$ and $\Delta$
is the comultiplication that defines the bialgebra
structure of~$\FU$. The corresponding solutions to \rf{RLL}
are expressed, respectively, in terms of the Gamma function
or its q--analogue (see \cite{J1,T1,F1,B2} for more details
in the latter case).

The aim of the present article is to develop a similar
algebraic construction of fundamental R--operators for models
whose underlying symmetry corresponds, in the sense of
Eq.~\rf{RLL0}, to the {\em quantum group}~$\Gq$.
More precisely, we introduce the quantum group $\Gqe$
with generators $a$, $b$, $c$, $d$, $\theta$, where
$\theta$ may be chosen to be the inverse to $b$ or~$c$.
It will be important to consider special {\em positive}
representations of $\Gqe$ which ensures that the operators
that we use are positive self--adjoint. These properties
are crucial for constructing fundamental R--operators
since we will need non--polynomial functions of generators
and their coproducts.

The article is organized as follows.
First, we discuss {\em Baxterization} of $\Gq$ and $\Gqe$,
presenting their defining relations in the form~\rf{RLL0}.
The two matrices, $g(\la)$ and $\hat{g}(\la)$, which play
the role of an \hbox{L--matrix} for~$\Gqe$, will be
our main objects of consideration.
Next, we show that, besides the standard comultiplication
$\Delta$, there is another algebra homomorphism
$\delta: \Gqe \to \Gqe^{\otimes 2}$. Further, we solve
Eq.~\rf{RLL} for $g(\la)$ and~$\hat{g}(\la)$.
The corresponding fundamental R--operators are given
(up to some twists) by powers of, respectively, $\Delta(bc)$ and
\hbox{$\delta(ad{-}qbc)$}. Next, we  show that the
\hbox{L--matrices} of the lattice Liouville model and the q--DST
model are nothing but $g(\la)$ and $\hat{g}(\la)$ with appropriately
chosen representations of generators. Using this observation,
we  construct the corresponding local lattice Hamiltonians.
Finally, we consider some reductions of $\Gqe$, including
the q--oscillator algebra $\Aq$ and the Weyl algebra~$\Wq$.
Following the same scheme, we introduce reductions
of $g(\la)$ and $\hat{g}(\la)$, and of $\Delta$ and~$\delta$,
and then construct the corresponding fundamental R--operators
by solving~Eq.~\rf{RLL}. We discuss relation of these
R--operators and of the corresponding local lattice Hamiltonians
to the Volterra model, the relativistic Toda model, and a
lattice regularization of the free field.

Let us remark that, although our construction based
on the use of the comultiplication structure yields
expressions for fundamental R--operators mainly as
powers of coproducts of some elements, it is
often useful to rewrite these expressions in terms
of the {\em quantum dilogarithm} function or, more
precisely, its self--dual form \cite{F2,F3} which is
suitable for dealing with the $|q|=1$ case.
A brief account on this function along with several
related statements which we use in the main text
are given in the Appendix.

\section{$\Gq$ and its Baxterization} 
\label{BGQ}

Let $q=e^{i\gamma}$, where $\gamma \in (0,\pi)$.
We will use the abbreviated notation $\Gq$ for the algebra of
regular functions on the quantum group, ${\it Fun}\bigl(\Gq\bigr)$
(see~\cite{V1,CP,KS}).
\begin{defn}\label{DEFGq}
$\Gq$ is a unital associative algebra with generators
$a$, $b$, $c$, $d$, and defining relations
\beq{defSL}
\begin{aligned}
{}  [a,d]=({q{-}q^{-1}}) \, b \, c \,,& \qquad [b,c] = 0\,,  \\
{}  a \, b = q \, b \, a \,,\quad
 a \, c = q \, c \, a \,,\quad &
 b \, d = q \, d \, b \,,\quad
 c \, d = q \, d \, c \,.
\end{aligned}
\end{equation}
$\Sq$ is the factor algebra of $\Gq$ over the ideal generated
by the relation \hbox{$ad-qbc = {\sf 1}$}.
\end{defn}
Following the \hbox{R--matrix} approach to quantum groups~\cite{FRT},
the generators of $\Gq$ can be assembled  into a matrix,
$g=\left(\begin{smallmatrix} a & b \\ c & d \end{smallmatrix}\right)$.
Then, by direct inspection of 16 quadratic exchange relations,
one can verify the following assertion (see, e.g.~\cite{CP,KS}).
\begin{lem}\label{ggSL}
The defining relations \rf{defSL} are equivalent to the
following relation
\beq{Rgg}
  R_{\1\2} \ g_{\1\3} \, g_{\2\3} =
  g_{\2\3} \, g_{\1\3} \ R_{\1\2} \,,
\end{equation}
where the auxiliary \hbox{R--matrix} is given by either of the
following matrices
\beq{Rconst}
 R^+ =  \left( \! \begin{array}{cccc}
  q \!\!\! & & & \\ [-2mm]
  & 1 &  & \\ [-1mm]
  &  q{-}q^{-1} \!\!\! &  1 & \\ [-2mm]
  & & & \!\! q
  \end{array} \! \right) \,, \qquad
  R^- = \bigl(R^+_{\2\1}\bigr)^{-1} =
 \left( \! \begin{array}{cccc}
  q^{-1} \!\!\! & & & \\ [-2mm]
  & 1 & q^{-1} {-} q & \\ [-1mm]
  &  &  1 & \\ [-2mm]
  & & & \!\!\!\! q^{-1}
  \end{array} \! \right) \,.
\end{equation}
\end{lem}

In what follows we will also need the following
spectral parameter dependent \hbox{R--matrices}
\begin{align}
\label{Rqosc2}
 \hat{R}(\la) & = \la R^{+} - \la^{-1} R^{-} =
 \left( \! \begin{array}{cccc}
  \varpi(q \la) \!\!\!\!\! & & & \\ [-2mm]
  & \varpi(\la) \!\! & \la^{-1} \varpi(q) & \\ [-1mm]
  & \la \, \varpi(q) \!\! & \varpi(\la) & \\ [-2mm]
  & & & \!\!\!\!\!  \varpi(q  \la)
  \end{array} \! \right) \,, \\
\label{Rqosc}
 R(\la) &= \lambda^{\frac12 \sigma_\3 \ot 1} \,
 \hat{R}(\la) \, \lambda^{-\frac12 \sigma_\3 \ot 1} =
 \left( \! \begin{array}{cccc}
  \varpi(q\la) \!\!\! & & & \\ [-2mm]
  & \varpi(\la) & \varpi(q) & \\ [-1mm]
  & \varpi(q) & \varpi(\la)  & \\ [-2mm]
  & & & \!\!\! \varpi(q\la)
  \end{array} \! \right)  \,,
\end{align}
where $\varpi(\la) \equiv \la-\la^{-1}$ and
$\sigma_{\!\3}=\left(\begin{smallmatrix} 1 & 0 \\ 0 & -1
   \end{smallmatrix}\right)$.

In the theory of quantum groups, the notion of
{\em Baxterization} was originally introduced by V.~Jones \cite{J2}
in the context of knot theory. It refers to the procedure of
constructing spectral parameter dependent solutions to the
Yang--Baxter equation out of solutions to the constant
(spectral parameter independent) Yang--Baxter equation.
An example is provided by the expression for $\hat{R}(\la)$
in terms of $R^\pm$ in formula~\rf{Rqosc2}. Analogously,
an L--matrix satisfying the RLL relation \rf{RLL0} can be regarded
as Baxterized if it is constructed from L--matrices that satisfy
the constant RLL relation. For instance, the \hbox{L--matrix} of
the XXZ model (see, e.g.,~\cite{F1}) has the form
\beq{xxz}
 L_{\rm\klein XXZ}(\la) = \la \, L_+ + \la^{-1} L_- \,,
 \end{equation}
where $L_\pm$ satisfy the constant RLL relation with
constant R--matrices given by~\rf{Rconst}.

In the theory of quantum integrable models it is crucial that
an R--matrix is spectral dependent (see Introduction), and so
the Baxterization procedure serves as quite a common technique
for constructing new solutions to the Yang--Baxter equation and
hence new integrable models. However, what concerns the
Baxterization of L--matrices, the vast majority of examples
occur in the cases where the symmetry $\FU$ is a quantum
algebra, typically the universal enveloping of a quantum
Lie algebra, like $\FU=\Uq$ for the XXZ model.

Quantum groups, in particular $\Gq$, are usually not considered
{}from the point of view of Baxterization of L--matrices.
In the present paper, we will try to fill this gap a bit.
Let us commence with observation that equation~\rf{Rgg}
can be Baxterized, albeit in a somewhat weaker sense
than it is usually meant. For this purpose, we assemble
the generators of $\Gq$ into two matrices:
\beq{gla0}
  \g(\la) = \left(\begin{matrix}  a & \la \, b \\
  \la^{-1} c & d \end{matrix}\right)\,,\qquad
  \hat{\g}(\la) = \left(\begin{matrix}
  \la^{-1} c &  \la^{-1} d \\
  \la \, a & \la \, b \end{matrix}\right) \,.
\end{equation}

\begin{propn}\label{RGG0}
Each of the following matrix relations
\begin{align}
\label{Rggl0}
  R_{\1\2}(\la) \ \g_{\1\3}(\la\mu) \, \g_{\2\3}(\mu) &=
  \g_{\2\3}(\mu) \, \g_{\1\3}(\la\mu) \ R_{\1\2}(\la) \,, \\
\label{Rhggl0}
  \hat{R}_{\1\2}(\la) \ \hat{\g}_{\1\3}(\la\mu) \,
  \hat{\g}_{\2\3}(\mu) &= \hat{\g}_{\2\3}(\mu) \,
 \hat{\g}_{\1\3}(\la\mu) \ \hat{R}_{\1\2}(\la) \,,
\end{align}
holds if and only if the elements $a$, $b$, $c$, $d$
satisfy the defining relations~\rf{defSL}.
\end{propn}

\proof
Matrices \rf{gla0} are related to each other and to
the matrix $g$ as follows
\beq{gla2b}
 \g(\la) = \lambda^{\frac12 \sigma_\3} \,
 g \ \lambda^{-\frac12 \sigma_\3} \,, \qquad
 \hat{\g}(\la) =
 \lambda^{-\frac12 \sigma_\3} \,
 \sigma_{\!\1} \, \g(\la) \, \lambda^{\frac12 \sigma_\3} \,,
\end{equation}
where $\sigma_{\!\1}=\left(\begin{smallmatrix} 0 & 1 \\ 1 & 0
   \end{smallmatrix}\right)$. Notice also that
\beq{Rsym}
 [R(\la) \,,\, \sigma_{\!\3} \ot 1 + 1 \ot \sigma_{\!\3}]=0 \,,
 \qquad [R(\la) \,,\, \sigma_{\!\1}\ot \sigma_{\!\1}] = 0 \,.
\end{equation}
Substituting the first of relations \rf{gla2b} into \rf{Rhggl0},
using the first of relations \rf{Rsym}, and taking into account
relation \rf{Rqosc} between $\hat{R}(\la)$ and $R(\la)$, it is
easy to see that \rf{Rggl0} is equivalent to the relation
$\hat{R}_{\1\2}(\la) g_{\1\3}  g_{\2\3} = g_{\2\3} g_{\1\3}
\hat{R}_{\1\2}(\la)$, which is nothing but a linear combination
of the $R{=}R_+$ and $R{=}R_-$ versions of Eq.~\rf{Rgg}.
Since $\lambda$ is arbitrary here, we conclude that \rf{Rggl0} is
equivalent to \rf{Rgg} and hence, by Lemma~\ref{ggSL}, to~\rf{defSL}.
Similarly, substituting the second relation in \rf{gla2b} into
\rf{Rhggl0}, using \rf{Rqosc} to replace $\hat{R}(\la)$ with $R(\la)$,
and then taking into account both relations \rf{Rsym}, it is easy
to see that \rf{Rhggl0} is equivalent to \rf{Rggl0}, and hence
to~\rf{defSL}.
\QED

The proof shows that the Baxterization in \rf{gla0} is not
a true one in the sense that it can be removed by the twist
transformations~\rf{gla2b}. Furthermore, for $\g(\la)$,
the transfer--matrix
\hbox{$T_\g(\la) = \mathrm{tr}_a
   \bigl(\g_{a,\SRN}(\la) \ldots \g_{a,\1}(\la)\bigr)$}
does not actually depend on $\la$ and thus it is not a
generating function for integrals of motion. However, the
corresponding transfer--matrix $T_{\hat{\g}}(\la)$ for
$\hat{\g}(\la)$ depends on $\la$ nontrivially, and
the operator coefficients $T_n$ in its
expansion, $T_{\hat{\g}}(\la)=\sum_n \la^n T_n$,
form a set of mutually commuting elements of~$\Gq^{\ot \SRN}$.
\pagebreak[4]

\section{$\Gqe$ and related lattice models} 
\label{GQE}

\subsection{Definition of $\Gqe$ and its Baxterization}
\label{EGQ}

Let us introduce the following extension of the quantum
group~$\Gq$.
\begin{defn}\label{DEFGqe}
$\Gqe$ is a unital associative algebra with generators
$a$, $b$, $c$, $d$, $\theta$, and defining relations
\rf{defSL} and
\beq{theta}
 a \, \theta = q^{-1} \, \theta \, a \,,\quad
 \theta \, d = q^{-1} \, d \, \theta \,, \quad
 [b,\theta] = 0 \,, \quad [\theta,c] = 0 \,.
\end{equation}
\end{defn}

\begin{lem}\label{cent2}
For a generic $q$, the center of $\Gqe$ is generated by the
following elements
\beq{eta}
 D_q \equiv a \, d - q \, b \, c \,, \qquad
 \eta'_q \equiv \theta \, b  \,, \qquad
 \eta''_q \equiv \theta \, c \,.
\end{equation}
\end{lem}
\proof
First, it is straightforward to check that $D_q$, $\eta'_q$,
and~$\eta''_q$ commute with the generators of~$\Gqe$.
Next, any central element $C$ of $\Gqe$
can be represented as a linear combination of monomials
\hbox{$a^{n_a} d^{n_d} b^{n_b} c^{n_c} \theta^{n_\theta}$},
where all $n$'s are non--negative integers. Commutativity of
$C$ with $b$, $c$, and $\theta$ implies that $n_a=n_d$.
Therefore, $C$ is equivalently represented as a
linear combination of monomials
\hbox{$D_q^{n} b^{m} c^{k} \theta^{l}$}. Commutativity of
$C$ with $a$ and $d$ implies that $m+k=l$.
Hence, using \rf{eta}, we conclude that $C$ is represented as
a linear combination of monomials
\hbox{$D_q^{n} (\eta'_q)^{m} (\eta''_q)^{k}$}.
\QED

\begin{lem}\label{Rgpgm}
The defining relations \rf{defSL} and~\rf{theta}
are equivalent to the following set of equations:
\begin{equation}
\label{Rgg2}
 R_{\1\2} \ g^\pm_{\1\3} \, g^\pm_{\2\3} =
  g^\pm_{\2\3} \, g^\pm_{\1\3} \ R_{\1\2} \,, \qquad
 R^+_{\1\2} \ g^+_{\1\3} \, g^-_{\2\3} =
  g^-_{\2\3} \, g^+_{\1\3} \ R^+_{\1\2} \,,
\end{equation}
where $g^+=\left(\begin{smallmatrix} \theta & 0 \\
  a\ & b \end{smallmatrix}\right)$ and
$g^-=\left(\begin{smallmatrix} c & d \\
  0 & 0 \end{smallmatrix}\right)$,
the auxiliary matrices $R^\pm$ are given by \rf{Rconst},
and $R$ in the first relation is either of them.
\end{lem}

\proof Direct inspection. \QED

Let us assemble the generators of $\Gqe$ into two matrices
\beq{gla}
  g(\la) = \left(\begin{matrix}  a & \la \, b \\
  \la \, \theta {+} \la^{-1} c & d \end{matrix}\right)\,,\qquad
  \hat{g}(\la) = \left(\begin{matrix}
  \la \, \theta {+} \la^{-1} c &  \la^{-1} d \\
  \la \, a & \la \, b \end{matrix}\right) \,.
\end{equation}

\begin{propn}\label{RGG}
Each of the following matrix relations
\begin{align}
\label{Rggl}
  R_{\1\2}(\la) \ g_{\1\3}(\la\mu) \, g_{\2\3}(\mu) &=
  g_{\2\3}(\mu) \, g_{\1\3}(\la\mu) \ R_{\1\2}(\la) \,, \\
\label{Rhggl}
  \hat{R}_{\1\2}(\la) \ \hat{g}_{\1\3}(\la\mu) \, \hat{g}_{\2\3}(\mu) &=
  \hat{g}_{\2\3}(\mu) \, \hat{g}_{\1\3}(\la\mu) \ \hat{R}_{\1\2}(\la) \,,
\end{align}
holds if and only if the elements $a$, $b$, $c$, $d$, $\theta$
satisfy the defining relations \rf{defSL} and~\rf{theta}.
\end{propn}

\proof
Notice that the second relation in \rf{gla2b} remains true
for $g(\la)$ and $\hat{g}(\la)$ given by~\rf{gla}. Therefore,
the same line of arguments as in the proof of Proposition~\ref{RGG0}
establishes equivalence of relations \rf{Rggl} and~\rf{Rhggl}.
Thus, it suffices to prove only~\rf{Rhggl}. For this aim
we observe that
\beq{gpgm}
  \hat{g}(\la) = \lambda \, g^+ + \lambda^{-1} \, g^- \,,
\end{equation}
where $g^\pm$ were defined in Lemma~\ref{Rgpgm}.
Substitute now \rf{Rqosc2} and \rf{gpgm} into \rf{Rhggl} and match
coefficients at different powers of $\la$ and~$\mu$.
It is not difficult to check that resulting matrix relations
are exactly those contained in~\rf{Rgg2}. (For the coefficient
at $\la^0 \mu^0$, we have to take into account the relation
$R^+ = P(R^-)^{-1}P$ along with the Hecke identity
\hbox{$R^+ - R^- = (q{-}q^{-1})P$}, where $P$ is
the permutation in ${\mathrm{Mat}}(2)^{\otimes 2}$, i.e.,
$P g^\pm_{\1\3} P = g^\pm_{\2\3}$.) Thus, relations
\rf{Rggl} and~\rf{Rhggl} are equivalent to \rf{Rgg2}, and
hence, by Lemma~\ref{Rgpgm}, to the defining relations of~$\Gqe$.
\QED

Unlike their $\Gq$ prototypes \rf{gla0}, matrices \rf{gla} are true
Baxterizations of $g^+$ and~$g^-$. Indeed, their q--determinants
(see, e.g., \cite{BT2}, Appendix~C) are
\beq{qdet}
 \qdet g(\la) = -  \qdet \hat{g}(\la) =
 D_q - q^{-1} \, \la^2 \, \eta_q' \,,
\end{equation}
which implies that the dependence of $g(\la)$ and
$\hat{g}(\la)$ on $\la$ cannot be removed by transformations
of the type~\rf{gla2b}.

Let us emphasize a close similarity between our L--matrices for
$\Gqe$ and those for $\Uq$. Indeed, $\hat{g}(\la)$ in \rf{gpgm}
and $L_{\rm\klein XXZ}$ in \rf{xxz} are constructed in the
same way from their constant counterparts and they satisfy
the RLL relations with the same auxiliary R--matrices. Such
a similarity seems quite natural in view of a duality
between $\Sq$ and~$\Uq$ (see \cite{CP,KS}). However, this similarity
is not absolute because the constant matrices $L_\pm$ in \rf{xxz}
are nondegenerate and generate the Borel subalgebras of $\Uq$,
whereas $g_-$ is degenerate and division of $\Gqe$ into
the subgroups generated by $g_\pm$ looks somewhat asymmetric.

\subsection{Standard and non--standard comultiplications
for $\Gqe$}\label{SNCP}

Recall that the linear homomorphism
$\Delta \colon \Gq \to \Gq^{\ot 2}$
defined on generators as follows
\beq{delSL}
  \begin{aligned}
{}  \Delta(a) = a \ot a + b \ot c \,,\quad & \quad
 \Delta(b) = a \ot b + b \ot d  \,,  \\
{}  \Delta(c) = c \ot a + d \ot c \,, \quad & \quad
 \Delta(d) = c \ot b + d \ot d  \,.
  \end{aligned}
\end{equation}
is a coassociative algebra homomorphism, i.e., its
homomorphism property
$
  \Delta(x\,y) = \Delta(x) \Delta(y)
$
is compatible with the defining relations \rf{defSL}, and
it satisfies the coassociativity property
\begin{equation}
\label{coass}
 (id \ot \Delta) \Delta(x) = (\Delta \ot id) \Delta(x) \,.
\end{equation}
The proof of these assertions is very simple in the
\hbox{R--matrix} approach due to an observation that
\rf{delSL} can be rewritten in the matrix form as follows:
\beq{gdel}
  (id \ot \Delta) \, g =  g_{\1\2} \, g_{\1\3}  \,.
\end{equation}
The fact that the Casimir element of $\Gq$ is a group--like element
w.r.t. the map $\Delta$, that is
\beq{DelCasSL}
  \Delta(D_q) = D_q \ot D_q \,,
\end{equation}
implies that the same map \rf{delSL} defines also a
coassociative algebra homomorphism for~$\Sq$.

$\Gq$ can be equipped with a bialgebra structure if, in addition
to the map $\Delta$, the linear homomorphism
$\epsilon \colon \Gq \to \BC$ is defined on generators as follows:
$
  \epsilon(g) =
 \left(\begin{smallmatrix} 1 & 0 \\ 0 & 1 \end{smallmatrix}\right)
$.
Then $\Delta$ and $\epsilon$ become {\em comultiplication} and
{\em counit} maps, respectively.

A natural question about the algebra $\Gqe$ is whether
we can introduce for it a comultiplication map, and, in
particular, whether we can extend the definition (\ref{delSL})
to~$\Gqe$. It appears that to define $\Delta(\theta)$
compatible with (\ref{defSL}) and (\ref{delSL}) in a purely
algebraic manner is not straightforward.
However, for our purposes it will be sufficient to define
$\Delta(\theta)$ for (special positive representations of)
real forms of certain factor algebras of~$\Gqe$.

\begin{defn}\label{SQER}
$\GqeR$ is a real form of $\Gqe$  equipped with an
anti--involution * defined on generators by
\beq{SLR}
  a^* = a \,, \quad b^*= b  \,, \quad
  c^* = c \,, \quad d^*= d \,, \quad \theta^* = \theta \,.
\end{equation}
$\GqeRp$ and $\GqeRpp$ are the factor algebras of $\GqeR$
over the ideals generated, respectively, by the relations
$\eta'_q = {\sf 1}$ and $\eta''_q = {\sf 1}$.
\end{defn}
Apparently, the algebras $\GqeRp$ and $\GqeRpp$
are isomorphic; the corresponding isomorphism map
$\iota$  is defined on generators as follows: $\iota(a)=a$,
$\iota(d)=d$, $\iota(b)=c$, $\iota(c)=b$, $\iota(\theta)=\theta$.

\begin{defn}\label{PIG}
Let $\CB$ be an algebra of linear operators acting on a Hilbert
space~$\CK$. Let $\FU$ stand for $\GqeRp$ or $\GqeRpp$.
An irreducible representation
$\pi \colon \FU \to \CB$ is called {\em positive} if the following
operators are {\em self--adjoint} and {\em strictly positive} on~$\CK$: \\
i) $\pi(x)$\, for\, $x=a,b,c,d,\theta,D_q$;\\
ii) $q^{\frac{1}{2}}\pi(a)\bigl(\pi(x)\bigr)^{-1}$\,
and\, $q^{\frac{1}{2}}\bigl(\pi(x)\bigr)^{-1}\pi(d)$\,
 for\, $x=b,c$.
\end{defn}

\begin{rem}\label{Rem1}
In Definition~\ref{PIG}, elements of $\FU$ are realized by
{\em unbounded} operators. Following \cite{W1,W2}, we will understand
the Weyl--type relations $x y = e^{i \gamma } y x$ in the defining
relations \rf{defSL} and \rf{theta} in the sense that, for a given
pair of positive self--adjoint operators $\pi(x)$ and $\pi(y)$, the
following unitary equivalence relations
$\pi(x)^{i t} \pi(y) \pi(x)^{-i t} = e^{- \gamma t} \pi(y)$ and
$\pi(y)^{i t} \pi(x) \pi(y)^{-i t} = e^{ \gamma t } \pi(x)$ hold
for all $t\in \mathbb R$ and admit analytic continuation to
complex values of~$t$.
\end{rem}

\begin{rem}\label{Rem2}
Condition $ii)$ in Definition~\ref{PIG} ensures that, for a pair of
generators $x$ and $y$ which satisfy the Weyl--type relation, the
sum $\pi(x)+\pi(y)$ is a positive self--adjoint operator. Indeed, let
$u$ and $v$ be positive self--adjoint operators satisfying relation
$u v = q^2 vu$. Then, in general, the sum \hbox{$u\,{+}\,v$} is
a symmetric but not necessarily self--adjoint operator~\cite{S1}.
If, following \cite{W1,W2}, we require that the operator
$q u^{-1} v$ is positive self--adjoint,
then property \rf{Sunit} of the quantum dilogarithm function
$S_\omega(t)$ (see Appendix~\ref{ApA}) implies that
$S_\omega(q u^{-1} v)$ is a unitary operator. In this case
Eq.~\rf{uvS1} shows that \hbox{$u\,{+}\,v$} is unitarily equivalent
to both $u$ and $v$ and hence is a positive self--adjoint operator.
Let us remark also that understanding relation $u v = q^2 vu$,
$q=i\gamma$ in the sense of Remark~\ref{Rem1} is equivalent to say
that $[\log u, \log v] = 2i\gamma$. Then, restricting our
consideration to the case $\gamma\in(0,\pi)$, again ensures
self--adjointness of \hbox{$u\,{+}\,v$}, by Proposition~A.2
in~\cite{S1}.
\end{rem}

An example of a positive representation of $\FU$ will be given below
in Section~3.5. Notice that $\pi(\eta'_q)$ and $\pi(\eta''_q)$
are also represented by positive self--adjoint operators.
Moreover, we have $\pi(\theta) = \bigl(\pi(b)\bigr)^{-1}$
for $\FU=\GqeRp$ and $\pi(\theta) = \bigl(\pi(c)\bigr)^{-1}$ for
$\FU=\GqeRpp$.

\begin{propn}\label{Bialg2}
Let $\CB$, $\CK$, and $\FU$ be as in Definition~\ref{PIG}
and let $\pi$ be a positive representation of~$\FU$.
Define the map $\Delta_\pi : \FU \to \CB^{\ot 2}$ as a linear
homomorphism such that:\\
i) $\Delta_\pi(x)=(\pi \ot \pi) \bigl(\Delta(x)\bigr)$ for
$x=a,b,c,d$ with $\Delta(x)$ given by~\rf{delSL}; \\
ii) $\Delta_\pi(\theta)= \bigl(\Delta_\pi(b)\bigr)^{-1}$
for $\FU=\GqeRp$ and
$\Delta_\pi(\theta)= \bigl(\Delta_\pi(c)\bigr)^{-1}$ for
$\FU=\GqeRpp$. \\
Then $\Delta_\pi$ is an {\em algebra homomorphism} and a
{\em \hbox{*--homomorphism}} w.r.t. the anti--involution~\rf{SLR}.
\end{propn}

\proof
The crucial property of $\Delta_\pi(x)$ for $x=a,b,c,d$ is that each
of these operators is of the form $u_x+v_x$, where $u_x$ and $v_x$
are positive self--adjoint operators satisfying the relation
$u_x v_x=q^2 v_x u_x$, e.g. $u_b=\pi(a){\ot}\pi(b)$ and
$v_b=\pi(b){\ot}\pi(d)$ for~$x=b$.
Furthermore, it is easy to check that $q u_x^{-1} v_x$ for $x=a,b,c,d$
are positive self--adjoint operators thanks to
condition $ii)$ in Definition~\ref{PIG}. According to Remark~\ref{Rem2},
these facts together imply that $\Delta_\pi(x)$ for $x=a,b,c,d$ are
also positive self--adjoint and hence {\em invertible} operators.
This, in particular, means that the inverse operators in the part
$ii)$ in the definition of~$\Delta_\pi$ are well defined.

Since $\Delta_\pi$ is a homomorphism, it suffices to verify its
properties for the generators. In particular, the
\hbox{*--homomorphism} property, which is
$\bigl(\Delta_\pi(x)\bigr)^* \equiv (* \ot *)\Delta_\pi(x)
 = \Delta(x^*)$ is obvious. The algebra homomorphism property of
$\Delta_\pi$ for $x=a,b,c,d$ is inherited {}from that of $\Delta$
for~$\Gq$. Finally, applying $\Delta_\pi{\ot}\Delta_\pi$
to \rf{theta} and multiplying the resulting relations with
$\Delta_\pi(b)$ (or $\Delta_\pi(c)$), we see that they are
equivalent to correct relations between $\Delta_\pi(b)$
(resp. $\Delta_\pi(c)$) and $\Delta_\pi(x)$ for $x=a,b,c,d$.
\QED

\begin{rem}
Using the $u_x+v_x$ form of $\Delta_\pi(x)$ along with Eq.~\rf{uvS1},
we can write an explicit expression for $\Delta_\pi(\theta)$.
For instance, in the case of $\FU=\GqeRp$ we have
\beq{delthb}
 \Delta_\pi(\theta) =
 S_\omega(\sw) \bigl(\pi(a)\,{\ot}\,\pi(b)\bigr)^{-1}
 \bigl(S_\omega(\sw)\bigr)^{-1} \,, \qquad
 \sw = \pi(b)\bigl(\pi(a)\bigr)^{-1} {\ot}
    \bigl(\pi(b)\bigr)^{-1}\pi(d) \,.
\end{equation}
\end{rem}

We introduced the map $\Delta_\pi$ by extending the standard
comultiplication~\rf{delSL} to~$\Gqe$. Now we will show that $\Gqe$
admits another ``comultiplication''
$\delta$ which is not related to~$\Delta$.

\begin{propn}\label{Bialg3}
The linear homomorphism $\delta \colon \Gqe \to \Gqe^{\ot 2}$
defined on generators as follows
\begin{eqnarray}
 &&
\begin{aligned}
{} & \delta(a) = a \ot \theta + b \ot a \,, \qquad
   \delta(\theta) = \theta \ot \theta \,,  \\
{} & \delta(c) = c \ot c  \,, \qquad
   \delta(b) = b \ot b \,, \qquad
 \delta(d) = c \ot d   \,,
\end{aligned}
 \label{delSL2}
\end{eqnarray}
is a {\em coassociative algebra homomorphism} and a
{\em \hbox{*--homomorphism}} w.r.t. the anti--involution~\rf{SLR}.
\end{propn}
\proof First, for the \hbox{*--homomorphism} property, it
suffices to notice that it obviously holds on generators.
Next, we notice that
\beq{delgh}
  (id \ot \delta) \, g^\pm =  g^\pm_{\1\2} \, g^\pm_{\1\3} \,,
\end{equation}
where $g^\pm$ were defined in Lemma~\ref{Rgpgm}.
This allows us to use the same approach as in the case
of~$\Gq$. Namely, the coassociativity property \rf{coass}
follows immediately if we apply
$\delta_{\2} \equiv (id \ot \delta \ot id)$ and
$\delta_{\3} \equiv (id \ot id \ot \delta)$ to~\rf{delgh}.
In order to prove compatibility of the homomorphism property
of $\delta$ with the defining relations \rf{defSL} and \rf{theta},
we recall that, by Lemma~\ref{Rgpgm}, these relations are
equivalent to relations~\rf{Rgg2}. Therefore, it suffices to
apply $\delta_{\3}$ to \rf{Rgg2}, use \rf{delgh}, and then to
verify the resulting R--matrix relations. The latter task
simply amounts to using \rf{Rgg2} twice, for instance:
$\delta_{\3} (R^+_{\1\2} g^+_{\1\3} g^-_{\2\3}) =
 R^+_{\1\2} g^+_{\1\3} g^+_{\1\4} g^-_{\2\3} g^-_{\2\4}=
 R^+_{\1\2} g^+_{\1\3} g^-_{\2\3} g^+_{\1\4}  g^-_{\2\4} =
 g^-_{\2\3} g^+_{\1\3} g^-_{\2\4} g^+_{\1\4}R^+_{\1\2} =
 g^-_{\2\3} g^-_{\2\4} g^+_{\1\3}  g^+_{\1\4}R^+_{\1\2}
 = \delta_{\3} ( g^-_{\2\3}  g^+_{\1\3}  R^+_{\1\2})$.
\QED

Notice that for $\delta$ there exists no counit $\epsilon$
because the bialgebra axiom
\hbox{$(id \ot \epsilon)\circ \delta = id\,$} cannot be
fulfilled as seen from the action of $\delta$ on~$d$.
Nevertheless, Proposition~\ref{Bialg3} justifies
referring to $\delta$  as a (non--standard) ``comultiplication''
for the sake of brevity.

An important difference of the non--standard ``comultiplication''
{}from $\Delta$ is that the generators $b$, $c$, and $\theta$
are group--like w.r.t.~$\delta$.
Therefore, so are the central elements~\rf{eta}:
\beq{deleta}
 \delta(\eta'_q) =\eta'_q \ot \eta'_q \,, \qquad
 \delta(\eta''_q) =\eta''_q \ot \eta''_q \,.
\end{equation}
On the other hand, the Casimir element $D_q$ is now not group--like.
Instead, we have
\beq{DelD2}
 \delta(D_q) = a c \ot \theta d + bc \ot D_q \,.
\end{equation}
Therefore, the relation $D_q={\sf 1}$ cannot be imposed as a
representation independent condition on generators.

Although both matrices $g(\la)$ and $\hat{g}(\la)$ define,
according to Proposition~\ref{RGG},
the same algebra $\Gqe$, the map $\delta$ is in a sense more
related to~$\hat{g}(\la)$. Indeed, formulae \rf{gpgm} and
\rf{delgh} have strong similarity with \rf{xxz} and
the formula $(id\ot\Delta)L_\pm=(L_\pm)_{\1\2} (L_\pm)_{\1\3}$,
which holds for the standard comultiplication of~$\Uq$.
We will see below that the construction of the fundamental
R--operator for $\hat{g}(\la)$ indeed requires invoking
the map~$\delta$, whereas the corresponding construction for
$g(\la)$ uses the map~$\Delta_\pi$.

\subsection{Fundamental R--operator for $g(\la)$}
\label{GQR1}

According to Proposition~\ref{RGG}, both matrices $g(\la)$
and $\hat{g}(\la)$ can serve as an \hbox{L--matrix} for
the algebra~$\FU=\Gqe$. Following the general scheme outlined
in Introduction, we have now to find their corresponding
fundamental R--operators, i.e., to solve Eq.~\rf{RLL}.
In this context, the following preliminary remark is in order.
In the case of $\FU=\Uq$, the L--matrices for the XXZ model and for
the sinh--Gordon model are related in essentially the same way as
$g(\la)$ and $\hat{g}(\la)$ (cf. the second relation
in Eq.~\rf{gla2b}) and, as a consequence, their fundamental
R--operators are also closely related~\cite{FT2,T1,BT2}. But in our
case there will be no such a relationship between the fundamental
R--operators for $g(\la)$ and~$\hat{g}(\la)$. To explain
this difference between our case and the $\Uq$ case,
let us formulate the following statement.

\begin{lem}\label{AUT}
Let $s$ be a constant invertible matrix. Suppose that matrices
$L(\la)$ and $\hat{L}(\la)= s \cdot L(\la)$
satisfy Eq.~\rf{RLL0} and define the same algebra~$\FU$.
If there exists an automorphism $\iota$ of $\FU$ such
that
\beq{aut}
 s \cdot L(\la) \cdot s = (id \ot \iota) \, L(\la) \,,
\end{equation}
then the fundamental R--operators corresponding to
$L(\la)$ and $\hat{L}(\la)$ are related as follows
\beq{autR}
 \SR(\la) = (\iota^{-1} \ot id) \, \hat{\SR}(\la) \,.
\end{equation}
\end{lem}

\proof
Consider Eq.~\rf{RLL} for $\hat{L}(\la)$, substitute all
$\hat{L}(\la)$ with $s \cdot L(\la)$, and use~\rf{aut}.
\QED

The structure of the L--matrices for the XXZ model and the
sinh--Gordon model is such that the automorphism $\iota$
does exist (for the generators of~$\Uq$ it reads:
$\iota(E)=F$, $\iota(F)=E$, $\iota(K)=K^{-1}$).
But for $g(\la)$ given by \rf{gla} and $s=\sigma_\1$,
matrix entries of the l.h.s. and the r.h.s. in \rf{aut}
have different functional dependence on~$\lambda$.
This means that there is no automorphism $\iota$
that would resolve \rf{aut} in our case and so
we have to solve Eq.~\rf{RLL} separately for $g(\la)$
and~$\hat{g}(\la)$.

Now we will solve Eq.~\rf{RLL} for~$g(\la)$.
For brevity of notations, we will write $x\,{\ot}\,y$
instead of $\pi(x)\,{\ot}\,\pi(y)$.

\begin{thm}\label{RG} 
Let $\CB$, $\mathcal K$, and $\FU$ be as in Definition~\ref{PIG}
and let $\pi$ be a positive representation of~$\FU$.
Let \hbox{$g(\la) \,{\in}\, {\mathrm{Mat}}(2){\ot}\CB$} be as
in~\rf{gla}. Then the operator $\SR(\la) \in \CB^{\ot 2}$ acting on
${\mathcal K}{\ot}{\mathcal K}$ and defined by the formula
\beq{Ral}
  \SR(\la) =  (c \ot b)^{-\frac{\al}{2} \log \la}\,
  \bigl( (a \ot b + b \ot d)(c \ot a + d \ot c) \bigr)^{\al \log \la}
   \, (c \ot b)^{-\frac{\al}{2} \log \la} \,,
\end{equation}
where~$ \al \equiv \fr{1}{\log q}=\fr{1}{i\gamma}$,
satisfies the equation
\beq{RLL1}
 \SR_{\2\3}(\la) \, g_{\1\2}(\la \mu) \, g_{\1\3} (\mu) =
 g_{\1\2} (\mu) \, g_{\1\3}(\la \mu) \, \SR_{\2\3}(\la) \,.
\end{equation}
If the tensor product $\pi \ot \pi$ is multiplicity free,
then \rf{Ral} is the unique solution of \rf{RLL1}
up to multiplication by a scalar factor.
\end{thm}

\proof
Matching coefficients at different powers of $\mu$, it is easy to
see that \rf{RLL1} is equivalent to the following set of equations:
\begin{align}
\label{Rcc}
{} [ \SR(\la) , (\theta \ot b) ] = 0 \,, & \qquad
{} [ \SR(\la) , (b \ot \theta) ] =0 \,, \\
\label{Rab}
{} \SR(\la) \, ( c \ot a + \la \, d \ot c)  &=
 (\la \, c \ot a + d \ot c) \, \SR(\la)\,, \\
{} \SR(\la) \, (a \ot b + \la\, b \ot d)  &=
 (\la\, a \ot b +  b \ot d) \, \SR(\la)\,, \\
 \label{Raa}
{} \SR(\la) \, (a \ot a + \la\, b \ot c)  &=
 (a \ot a + \la^{-1} b \ot c) \, \SR(\la)\,, \\
\label{Rdd}
{} \SR(\la) \, (d \ot d + \la^{-1} c \ot b )  &=
 (d \ot d + \la\, c \ot b ) \, \SR(\la)\,, \\
{} \SR(\la) \, (\la\, \theta \ot a + d \ot \theta )  &=
 (\theta \ot a + \la\, d \ot \theta) \, \SR(\la) \,.
 \label{Rac1}
\end{align}
It is now easy to recognize in \rf{Rab}--\rf{Rdd} a structure
related to the comultiplication~$\Delta$ (cf.~\rf{delSL}).
To make this structure more transparent, we introduce
$\tilde{\SR}(\la) = (c \ot b)^{\frac{\al}{2} \log \la} \,
 \SR(\la) \, (c \ot b)^{\frac{\al}{2} \log \la}$.
Then equations \rf{Rcc}--\rf{Rac1} acquire the following form:
\begin{align}
\label{Rcc2}
{} [ \tilde{\SR}(\la) , b \ot \theta ] &=
 [ \tilde{\SR}(\la) , \theta \ot b] =0 \,, \\
 \label{Rbc}
 \tilde{\SR}(\la) \, \Delta_\pi(b) =
 \Delta_\pi(b) \, \tilde{\SR}(\la) \,,  \quad & \qquad
 \tilde{\SR}(\la) \, \Delta_\pi(c) =
 \Delta_\pi(c) \, \tilde{\SR}(\la) \,, \\
 \label{Rad}
 \tilde{\SR}(\la) \, \Delta_\pi(a) =
 \la^{-2} \Delta_\pi(a) \, \tilde{\SR}(\la) \,, \quad & \qquad
 \tilde{\SR}(\la) \, \Delta_\pi(d) =
 \la^2 \, \Delta_\pi(d) \, \tilde{\SR}(\la) \,, \\
\label{Rac2}
 \tilde{\SR}(\la) \, (\la \, \theta \ot a
    + \la^{-1} d \ot \theta )  &=
 (\la^{-1} \theta \ot a  + \la \, d \ot \theta ) \,
    \tilde{\SR}(\la) \,.
\end{align}
where $\Delta_\pi$ is the algebra homomorphism introduced
in Proposition~\ref{Bialg2}. Next, observing that
\begin{align}
\label{adc}
 [ \Delta_\pi(a), b \ot \theta]=0 \,, \quad & \qquad
  \Delta_\pi(b) \, (b \ot \theta) =
  q \,(b \ot \theta) \, \Delta_\pi(b) \,, \\
\label{bcc}
[ \Delta_\pi(d), b \ot \theta]=0
  \,,  \quad & \qquad
  \Delta_\pi(c) \, (b \ot \theta) =
  q^{-1} (b \ot \theta) \, \Delta_\pi(c) \,,
\end{align}
are consequences of \rf{defSL}, \rf{theta}, and \rf{delSL},
we infer that Eqs.~\rf{Rcc2}--\rf{Rbc} are satisfied if
$\tilde{\SR}(\la)$ is taken to be a function of
$\Delta_\pi(ad)$ and $\Delta_\pi(bc)$. Furthermore, due to
Eq.~\rf{DelCasSL} we have
$\Delta_\pi(ad)= q\Delta_\pi(bc){+} D_q\,{\ot}\,D_q$,
where the last term is a multiple of the unit operator.
This implies that we can take $\tilde{\SR}(\la)$ to be a function of
$\Delta_\pi(bc)$ only. Then Eqs.~\rf{Rcc2}--\rf{Rad} are solved easily:
\beq{Rl1}
 \tilde{\SR}(\la) =  \bigl(\Delta_\pi(bc)\bigr)^{\al \log \la}
 \,, \qquad  \al = \fr{1}{\log q} \,.
\end{equation}
It remains to verify \rf{Rac2}. For this aim we notice that,
since $\Delta_\pi(b)$ is invertible, Eq.~\rf{Rac2} is equivalent
to the relation
\beq{rx}
 \tilde{\SR}(\la) \, \SX(\la) =
 \SX(\la^{-1}) \, \tilde{\SR}(\la) \,,
\end{equation}
where we denoted
$\SX(\lambda) \equiv  \Delta_\pi(b)
    (\la \, \theta \ot a + \la^{-1} d \ot \theta )$.
Now, using \rf{delSL}, we find
\beq{X1}
 \SX(\la) = q^{-1} \la (\theta \ot b) \, \Delta_\pi(a) +
 q\, \la^{-1}  (b \ot \theta) \, \Delta_\pi(d)
  + \la \, \eta'_q \ot D_q  + \la^{-1} \, D_q \ot \eta'_q \,.
\end{equation}
The sum of the last two terms here obviously satisfies~\rf{rx}.
The first two terms satisfy \rf{rx} as a consequence of
relations \rf{Rcc2} and~\rf{Rad}. Thus, Eq.~\rf{Rac2} is proven
and we have shown that \rf{Ral} indeed solves Eq.~\rf{RLL1}.

Let us prove the uniqueness of~$\SR(\la)$. Notice that
$\tilde{\SR}(\la)$ is an invertible operator due to the
properties of~$\pi$. Suppose that there exists another solution,
$\tilde{\SR}'(\la)$, to Eqs.~\rf{Rcc2}--\rf{Rad}. Then it
follows from \rf{Rbc}--\rf{Rad} that
$\SF(\la) \equiv (\tilde{\SR}(\la))^{-1}\tilde{\SR}'(\la)$
commutes with $\Delta_\pi(x)$ for all $x \in \FU$.
Under the assumption that $\pi\,{\ot}\,\pi$ is multiplicity free,
we invoke Lemma~\ref{cent2} and infer that $\SF(\la)$ can
be a function only of~$\Delta_\pi(\eta''_q)$.
But it follows from \rf{adc}--\rf{bcc} that $\Delta_\pi(\eta''_q)$
does not commute with $b\,{\ot}\,\theta$.
Thus, $\SF(\la)$ satisfying \rf{Rcc2} cannot depend non--trivially
on $\Delta_\pi(\eta''_q)$ and therefore it must be just
a scalar function.
\QED

\begin{rem}
The positivity property of the representation~$\pi$ is crucial
for the assertion that \rf{Rl1} solves Eqs.~\rf{Rbc}--\rf{Rad}.
Indeed, it ensures that $x=\Delta_\pi(bc)$ and $y=\Delta_\pi(z)$
for $z=a,b,c,d$ are positive self--adjoint operators
(cf. Remark~\ref{Rem2}) and therefore  \rf{Rl1} solves
Eqs.~\rf{Rbc}--\rf{Rad} in the sense clarified in Remark~\ref{Rem1}.
Notice also that on the same ground we have
$\bigl(\Delta_\pi(b)\Delta_\pi(c)\bigr)^t=\bigl(\Delta_\pi(bc)\bigr)^t$.
\end{rem}

\begin{rem}
For lattice integrable models, the function that most commonly
appears in solutions for fundamental R--operators is the
{\em quantum dilogarithm} (see Appendix~\ref{ApA}).
Lemma~\ref{LUL} (see the same Appendix) allows us to rewrite
our solution \rf{Ral} in a form involving quantum dilogarithms:
\beq{Raldil}
\begin{aligned}
  \SR(\la) &= \frac{ S_\omega(\la^{-1} \sw)}{S_\omega(\la\, \sw)} \,
 (a \ot a)^{\alpha \log\la} \,
 \frac{S_\omega(\la^{-1} \tilde{\sw})}{S_\omega(\la\, \tilde{\sw})} \\
 {}&= \frac{S_\omega(\la^{-1} \sw^{-1})}{ S_\omega(\la\, \sw^{-1})} \,
 (d \ot d)^{\alpha \log\la} \,
 \frac{S_\omega(\la^{-1} \tilde{\sw}^{-1})}%
    {S_\omega(\la\, \tilde{\sw}^{-1})} \,,
\end{aligned}
\end{equation}
where $\sw=b a^{-1} {\ot} b^{-1}d$ and
$\tilde{\sw}=d c^{-1} {\ot} a^{-1}c$.
\end{rem}

Fundamental R--operator \rf{Ral} is regular in the sense
of Eq.~\rf{reg} and has the following properties:
\beq{unit}
 \bigl(\SR(\bar{\la})\bigr)^* =\SR(\la^{-1}) = \SR^{-1}(\la)  \,.
\end{equation}
Application of formula \rf{Hloc} to \rf{Ral} yields the
following lattice Hamiltonian density:
\beq{Hloc1}
 {\gamma} \, H_{n,n \p1} =
 \log \bigl( (a_{n \p1} b_{n} + b_{n \p1} d_{n})
 (c_{n \p1} a_{n} + d_{n \p1} c_{n}) \bigr) - \log(b_n b_{n \p1}) \,.
\end{equation}
Definition~\ref{PIG} along with Remark~\ref{Rem2} ensure that
the arguments of the logarithms here are products of commuting
positive self--adjoint operators.

\subsection{Fundamental R--operator for $\hat{g}(\la)$}

Now we will solve Eq.~\rf{RLL} for $\hat{g}(\la)$.
For brevity of notations, we will write $x\,{\ot}\,y$
instead of $\pi(x)\,{\ot}\,\pi(y)$ and $\delta$ instead
of~$(\pi\,{\ot}\,\pi) \circ \delta$.

\begin{thm}\label{RGh} 
Let $\CB$, $\mathcal K$, and $\FU$ be as in Definition~\ref{PIG}
and let $\pi$ be a positive representation of~$\FU$.
Let \hbox{$\hat{g}(\la) \,{\in}\, {\mathrm{Mat}}(2){\ot}\CB$} be as
in~\rf{gla}. Then the operator $\hat{\SR}(\la) \in \CB^{\ot 2}$
acting on ${\mathcal K}{\ot}{\mathcal K}$ and defined by the formula
\beq{Ralh}
  \hat{\SR}(\la) =
  \bigl( a c \ot \theta d  +
   b c \ot a d  - q bc \ot bc \bigr)^{\al \log \la} \,,
\end{equation}
where~$ \al \equiv \fr{1}{\log q}=\fr{1}{i\gamma}$,
satisfies the equation
\beq{RLL2}
 \hat{\SR}_{\2\3}(\la) \, \hat{g}_{\1\2}(\la \mu) \,
  \hat{g}_{\1\3} (\mu) = \hat{g}_{\1\2} (\mu) \,
  \hat{g}_{\1\3}(\la \mu) \, \hat{\SR}_{\2\3}(\la) \,.
\end{equation}
If the tensor product $\pi \ot \pi$ is multiplicity free, then
\rf{Ralh} is the unique solution of \rf{RLL2} up to
multiplication by a scalar factor.
\end{thm}

\proof
Substituting the Baxterized form \rf{gpgm} of $\hat{g}(\la)$
into \rf{RLL2} and matching coefficients at different powers of
$\mu$, we see that \rf{RLL2} is equivalent to the following
set of matrix equations:
\begin{align}
\label{RLL2gg}
 \hat{\SR}_{\2\3}(\la) \, \hat{g}^\pm_{\1\2} \,
  \hat{g}^\pm_{\1\3} &= \hat{g}^\pm_{\1\2}  \,
  \hat{g}^\pm_{\1\3} \, \hat{\SR}_{\2\3}(\la) \,,\\
\label{RLL2ggg}
 \hat{\SR}_{\2\3}(\la) \, (\la \, \hat{g}^+_{\1\2} \,
  \hat{g}^-_{\1\3} + \la^{-1} \hat{g}^-_{\1\2} \, \hat{g}^+_{\1\3} )
 &= (\la^{-1}  \hat{g}^+_{\1\2} \, \hat{g}^-_{\1\3} +
 \la \, \hat{g}^-_{\1\2} \, \hat{g}^+_{\1\3} ) \,
 \hat{\SR}_{\2\3}(\la) \,.
\end{align}
Comparing \rf{RLL2gg} with \rf{delgh}, we conclude that
\beq{SRd}
 [\, \hat{\SR}(\la) \,,\, \delta(x) \,] = 0
\end{equation}
for all generators and hence for all~$x\in \FU$. This suggests
to seek $\hat{\SR}(\la)$ as a function of~$\delta(D_q)$.

Matrix equation \rf{RLL2ggg} is equivalent to the
following set of equations:
\begin{align}
\label{Rabd}
{} \hat{\SR}(\la) \, (a \ot d ) =
 \la^{-2} (a \ot d ) \, \hat{\SR}(\la) \,, \ & \ \quad
 \hat{\SR}(\la) \, (a \ot c) =
 \la^{-2} (a \ot c) \, \hat{\SR}(\la) \,,  \\
 \label{Re1}
{} \hat{\SR}(\la) \, (\la\, \theta \ot d + \la^{-1} d \ot b) &=
 (\la^{-1} \theta \ot d + \la\, d \ot b) \, \hat{\SR}(\la)\,,\\
\label{Re2}
{} \hat{\SR}(\la) \, (\la^{-1} d \ot a + \la^{-1} c \ot \theta +
 \la\, \theta \ot c) &=
 (\la\, d \ot a + \la\, c \ot \theta +
 \la^{-1} \theta \ot c) \, \hat{\SR}(\la) .
\end{align}
Noticing that
\beq{Dabd}
 \delta(D_q) \, (a\ot d) =  q^{-2} (a\ot d) \, \delta(D_q) \,, \qquad
 \delta(D_q) \, (a\ot c) =  q^{-2} (a\ot c) \, \delta(D_q) \,,
\end{equation}
we infer that a solution to \rf{Rabd} is given by
\beq{Rsol}
 \hat{\SR}(\lambda) =  \bigl( \delta(D_q) \bigr)^{\al \log \la}
 \,, \qquad  \al \equiv \fr{1}{\log q} \,.
\end{equation}

\begin{lem}\label{RBD}
 $\hat{\SR}(\lambda)$ given by \rf{Rsol} satisfies relation~\rf{Re1}.
\end{lem}
\noindent
The proof is given in Appendix~B.
It remains to prove that $\hat{\SR}(\lambda)$ satisfies Eq.~\rf{Re2}.
For this aim, we notice that since $\delta(b)$ is represented by an
invertible element, Eq.~\rf{Re2} is equivalent to
the relation
\beq{rxx}
 \hat{\SR}(\la) \, \SX(\la) =
 \SX(\la^{-1}) \, \hat{\SR}(\la) \,,
\end{equation}
where we denoted
$\SX(\lambda) \equiv  q^{-1} \delta(b) \,
(\la^{-1} d \ot a + \la^{-1} c \ot \theta + \la\, \theta \ot c) $
Now we observe that
\beq{X2}
 \SX(\la) =
 (\la \, \theta \ot d + \la^{-1} d \ot b) \, \delta(a)  -
 q \, \la \, \delta(\theta)\, (a \ot d)  -
  \la \, D_q \ot \eta'_q - \la^{-1} \eta'_q \ot D_q  \,.
\end{equation}
The sum of the last two terms here obviously satisfies~\rf{rxx}.
The first two terms satisfy \rf{rxx} as a consequence of relations
\rf{Rabd} and~\rf{Re1}. Thus, \rf{Re2} is proven and we
have shown that \rf{Rsol} indeed solves Eq.~\rf{RLL2}.

Eq.~\rf{SRd} implies that $\hat{\SR}(\la)$ is essentially unique.
Indeed, under the assumption that $\pi \ot \pi$ is multiplicity free,
we invoke Lemma~\ref{cent2} and infer that $\hat{\SR}(\la)$ can
be a function only of $\delta(D_q)$ and~$\Delta_\pi(\eta''_q)$.
Furthermore, Eq.~\rf{deleta} implies that
$(\pi{\ot}\pi)\delta(\eta''_q)$ is just a multiple of unity,
so $\hat{\SR}(\la)$ must be a function of $\delta(D_q)$ only.
Finally, it is clear that such a function satisfying \rf{Dabd}
is given by \rf{Rsol} uniquely up to a scalar factor.
\QED

\begin{rem}
Lemma~\ref{LUL} allows us to rewrite our solution \rf{Ralh} in
terms of quantum dilogarithms:
\beq{Ralhdil}
  \hat{\SR}(\la) = \bigl( bc \ot D_q \bigr)^{\al \log \la} \
 \frac{ S_\omega( {\mathsf r} )}%
    {S_\omega\bigl(\la^2 {\mathsf r} \bigr)} \,,\qquad
 {\mathsf r} = \bigl(D_q\bigr)^{-1} \, b^{-1} a \ot \theta d   \,.
\end{equation}
\end{rem}

Fundamental R--operator \rf{Ralh} is regular and has the
properties~\rf{unit}. Application of formula \rf{Hloc} yields
the following lattice Hamiltonian density:
\beq{Hloc2}
 \gamma \, \hat{H}_{n,n \p1} =
 \log \bigl( a_{n \p1} c_{n \p1} \theta_n d_n  +
   b_{n \p1} c_{n \p1} (D_q)_n \bigr) \,.
\end{equation}
Definition~\ref{PIG} along with Remark~\ref{Rem2} ensure that the
argument of the logarithm here is a positive self--adjoint operator.

\subsection{Lattice Liouville model}

A one--parameter family $\pi_{\kappa}$ of positive representations
of $\FU=\GqeRp$ or $\FU=\GqeRpp$ on the Hilbert space $\CK = L^2(\BR)$
can be constructed as follows (it is closely related to one of the
representations of $\SqR$ listed in~\cite{S1}):
\beq{EFK1}
\begin{aligned}
{} \pi_{\kappa}(a)  =  e^{\frac{\beta}{8}\, \Pi} \,
    \bigl( \kappa^{-1} +\kappa  \, e^{-\beta \Phi} \bigr) \,
    e^{\frac{\beta}{8}\, \Pi} \,, \qquad &
 \pi_{\kappa}(d) = \kappa^{-1} e^{- \frac{\beta}{4} \, \Pi} \,, \\
{} \pi_{\kappa}(b) = \pi_{\kappa}(c) = e^{-\frac{\beta}{2} \, \Phi}
 \,, \qquad & \pi_{\kappa}(\theta) = e^{\frac{\beta}{2} \, \Phi} \,.
\end{aligned}
\end{equation}
Here $\kappa>0$, $\beta \equiv \sqrt{8\gamma} =2\omega\sqrt{2\pi} >0$,
and $\Pi$ and  $\Phi$ are self--adjoint operators
on $L^2(\BR)$ which satisfy \hbox{$[\Pi,\Phi]=-i$}.
Since elements of $\FU$ are realized by unbounded operators
on $L^2(\BR)$, it is necessary to consider suitable subspaces
$\CT_\kappa\subset L^2(\BR)$ of test--functions on which all operators
$\pi_\kappa(x)$,  $x\in\FU$ are well defined. Similar consideration
was done for $\UqR$ in~\cite{PT,BT1}. We will provide analogous
analytic details for $\pi_{\kappa}$ elsewhere.

Let us now introduce the following \hbox{L--matrix}:
$ L^{\kl}(\la) = \kappa\, \pi_\kappa \bigl( g(\la) \bigr)$.
In order to construct the corresponding lattice model
we assign a copy of this matrix to
each site of the lattice, i.e. for \hbox{$n=1,\ldots,\SRN$}
we have
\beq{LL}
 L_n^{\kl}(\la)   =
 \left( \begin{array}{cc}
  e^{\frac{\beta}{8}\, \Pi_n} \,
    \bigl( 1 +\kappa^2  \, e^{-\beta \Phi_n} \bigr) \,
    e^{\frac{\beta}{8}\, \Pi_n} &
  \kappa\, \la \, e^{ - \frac{\beta}{2} \Phi_n }  \\
  \kappa\, \bigl(\la \, e^{\frac{\beta}{2} \Phi_n} +
  \la^{-1}  e^{-\frac{\beta}{2} \Phi_n} \bigr) &
   e^{- \frac{\beta}{4} \, \Pi_n}
 \end{array} \right) \,,
\end{equation}
where $\Phi_n$ and $\Pi_n$ act non--trivially only on the $n$--th tensor
factor in the Hilbert space \hbox{$\CH=\bigl(L^2(\BR))^{\ot \SRN}$} and
therefore satisfy the relation \hbox{$[\Pi_n,\Phi_m]=-i \delta_{nm}$}.

In the pioneering work \cite{FT3}, a close analogue of \rf{LL} was
constructed as a special limit of the \hbox{L--matrix} for the
sine--Gordon model and put forward as an \hbox{L--matrix} describing
a lattice version of the Liouville model with $\Phi_n$ and $\Pi_n$ being
discrete counterparts of the field and its conjugate momentum variables.
In its present form, L--matrix \rf{LL} was obtained in \cite{BT2}
by analogous limit applied to the sinh--Gordon model.

The {\em continuum limit} of a classical lattice integrable model is
usually constructed as the limit of vanishing lattice spacing
($\SRN\to\infty$, $\kappa\to 0$ with $\kappa\,\SRN$ kept fixed)
combined with the standard recipe \cite{FST} of replacement of
lattice canonical variables by their continuum counterparts:
\beq{cont}
 \Pi_n \to \kappa \, \Pi(x) \,, \quad
 \Phi_n \to \Phi(x) \,, \quad
 x = n \, \kappa  \,,
\end{equation}
which leads to the canonical Poisson brackets,
$\{\Pi(x),\Phi(y)\}=\delta(x-y)$.
In this classical continuum limit we have
$L_n^{\kl}(\la) =
 \bigl(\begin{smallmatrix} 1 & 0\\ 0 & 1
    \end{smallmatrix} \bigr)
 + \kappa \, \bigl(U_+(\la)+U_-(\la)\bigr) + O(\kappa^2)$,
where ($\partial_\pm \equiv \partial_t \pm \partial_x$)
\beq{UccL}
  U_+(\la) =   \left( \begin{array}{cc}
  \fr{\beta}{8} \, \partial_+ \Phi  &
  \la \, e^{ - \frac{\beta}{2} \Phi }  \\
 \la \, e^{\frac{\beta}{2} \Phi}   &
  - \fr{\beta}{8} \, \partial_+ \Phi
 \end{array} \right) \,, \quad
 U_-(\la) =   \left( \begin{array}{cc}
  \fr{\beta}{8} \, \partial_- \Phi  &  0  \\
  \la^{-1}  e^{-\frac{\beta}{2} \Phi}  &
  - \fr{\beta}{8} \, \partial_- \Phi
 \end{array} \right) \,.
\end{equation}
These matrices satisfy the zero curvature equation,
$\partial_-U_+(\la) + \partial_+U_-(\la) = 2[U_+(\la),U_-(\la)]$,
provided that $\Phi$ satisfies the equation of motion of the Liouville
field: $\square\, \Phi = \frac{8}{\beta} e^{ - \beta \Phi }$.
On this ground it was suggested in \cite{FT3} that \rf{LL} corresponds
to a lattice version of the Liouville model. However, a direct
verification of this claim, i.e., construction of a lattice Hamiltonian
that $i)$ commutes with the transfer matrix for \rf{LL}, and $ii)$
turns in the continuum limit into the Hamiltonian of the continuum
Liouville model, has been missing until now although some partial results
have been obtained. In particular, it was shown in \cite{BT2} that
applying to the  Hamiltonian of the lattice sinh--Gordon model first
the special limit procedure described in \cite{FT3} and then taking
the continuum limit, we indeed obtain the Hamiltonian of the continuum
Liouville model. Another computation~\cite{S3} demonstrated that,
unlike for the sinh--Gordon model, the factorization method \cite{IK}
of constructing integrals of motion applied to \rf{LL} yields a lattice
analogue only of the chiral combination \hbox{$(H+P)$} of the Liouville
Hamiltonian and momentum operator.

Results of Section~3.3 imply that
\beq{HLq}
 H^{\kl}_{n,n \p1}=
 (\pi_\kappa \ot \pi_\kappa) H_{n,n \p1} \,,
\end{equation}
where $H_{n,n \p1}$ is given by \rf{Hloc1},
is a quantum nearest--neighbour lattice Hamiltonian corresponding
to L--matrix~\rf{LL}. In order to show that \rf{HLq} is a lattice
analogue of the Hamiltonian for the continuous Liouville model
we consider first its classical limit where $\Phi_n$ and $\Pi_n$
become canonical variables on the phase space equipped with the
Poisson bracket \hbox{$\{\Pi_n,\Phi_m\}=\delta_{nm}$}.
A direct computation using \rf{EFK1} yields
(up to an additive constant)
\begin{align}\label{Hcl}
 H^{\rm \kl, cl}_{n,n \p1}
 =&\ \fr{1}{\gamma}\log \Bigl(
 \fr{1}{2}\cosh \fr{\beta}{4} \, (\Pi_n + \Pi_{n \p1})
 + \fr{1}{2} \,
 \cosh \fr{\beta}{2} \, (\Phi_n - \Phi_{n \p1}) \\
\nonumber
 &+ \fr{\kappa^2}{2} \, e^{-\frac{\beta}{2}(\Phi_n + \Phi_{n \p1})} \,
 \bigl( 1 + e^{\frac{\beta}{4}(\Pi_n + \Pi_{n \p1})} \,
 \cosh \fr{\beta}{2} ( \Phi_n - \Phi_{n \p1} ) \bigr) \\
 \nn & + \fr{\kappa^4}{4} \,
 e^{\frac{\beta}{4}(\Pi_n + \Pi_{n \p1})}
 e^{-\beta(\Phi_n + \Phi_{n \p1})} \Bigr)\,.
\end{align}
Let us remark that the difference between \rf{Hcl} and the analogous
expression obtained by a ``naive'' limit in \cite{BT2} is only in the
last term.
Taking now the continuum limit of \rf{Hcl} according to \rf{cont},
we obtain (again up to an additive constant)
\begin{equation}
 \lim_{\kappa\to 0} \sum_n \fr{1}{\kappa}
 H^{\rm \kl, cl}_{n,n \p1} =
 \int dx \, \bigl( \fr{1}{2} \, \Pi^2 + \fr{1}{2} \,
 (\partial_x\Phi)^2
 + \fr{1}{\gamma} \, e^{-\beta \,\Phi} \bigr)  \,,
\end{equation}
which is the Hamiltonian of the classical continuum Liouville
model.

\section{Reductions of $\Gqe$ and related lattice models} 
\label{RGQE}

Defining relations of $\Gqe$ admit the following reductions:
$i)$ $\theta=0$; $ii)$ $b=c$; $iii)$ $b=0$, and $iv)$ $c=0$.
Below we will consider the problem of constructing the
fundamental R--matrices for the corresponding reductions
of matrices $g(\la)$ and $\hat{g}(\la)$ in each of these
cases.

\subsection{$\theta=0$} 
\label{TZ}

For $\theta=0$, matrix $g(\la)$ reduces back to $\g(\la)$ given
by~\rf{gla0}. In this case we take $\pi$ to be a positive
representation of~$\Gq$ (modification of Definition~\ref{PIG}
is obvious).
As we discussed at the end of Section~\ref{BGQ}, dependence on the
spectral parameter of the auxiliary transfer--matrix for $\g(\la)$
can be removed with the help of the twist transformation~\rf{gla2b}.
However, Eq.~\rf{RLL1} for the corresponding fundamental R--matrix
cannot be transformed by a similar means to a spectral parameter
independent form.

Let $\SR_\0(\la)$ denote a solution to \rf{RLL1} where $g(\la)$ is
replaced with $\g(\la)$.
Introduce $\tilde{\SR}_\0(\la) = (c \ot c )^{\frac{\al}{2} \log \la} \,
 \SR_\0(\la) \, (c \ot c)^{\frac{\al}{2} \log \la}$.
Evidently, $\SR_\0(\la)$ must satisfy only Eqs.~\rf{Rab}--\rf{Rdd},
and $\tilde{\SR}_\0(\la)$ must satisfy only Eqs.~\rf{Rbc}--\rf{Rad}.
For the latter we have a one--parameter family of solutions:
\beq{R0la}
\tilde{\SR}_\0(\la;\beta) =
 \bigl(\Delta_\pi(b)\bigr)^{(\al-\beta) \log \la}
 \bigl(\Delta_\pi(c)\bigr)^{(\al+\beta) \log \la} \,, \qquad
 \al = \fr{1}{\log q} \,.
\end{equation}

\begin{rem}
The reason why the proof of essential uniqueness given for
$\tilde{\SR}(\la)$ in Section~\ref{GQR1} does not apply to the case of
$\tilde{\SR}_\0(\la)$ despite that, by Lemma~\ref{cent2}, the center
of $\Gq$ is generated only by the quadratic Casimir element~$D_q$,
is that the ratio of two solutions,
$\SF(\la) = \bigl( \tilde{\SR}_\0(\la;\beta_\1) \bigr)^{-1}
 \tilde{\SR}_\0(\la;\beta_\2)  =
 \bigl(\Delta_\pi(bc^{-1})\bigr)^{(\beta_\1-\beta_\2) \log \la}$
depends non--trivially on the {\em non--polynomial}
Casimir element, $bc^{-1}$, which can formally
be written as $\lim_{\theta\to 0} \eta'/\eta''$.
\end{rem}

Next we consider the $\theta=0$ counterpart of matrix
$\hat{g}(\la)$ which is $\hat{\g}(\la)$ given by~\rf{gla0}.
We again take $\pi$ to be a positive representation of~$\Gq$.

Let $\hat{\SR}_\0(\la)$ denote a solution to \rf{RLL2} where
$\hat{g}(\la)$ is replaced with~$\hat{\g}(\la)$. Apparently,
$\hat{\SR}_\0(\la)$ must satisfy \rf{RLL2gg} and as a consequence
it is a function of $\delta(D_q)$ only (the Casimir element
$bc^{-1}$ is group--like w.r.t. $\delta$ and hence is represented
by a multiple of the unity). Further, $\hat{\SR}_\0(\la)$
must satisfy Eq.~\rf{Rabd} and the relations that replace
Eqs.~\rf{Re1}--\rf{Re2}, namely
$\hat{\SR}_\0(\la)\, (d {\ot} x ) =
 \la^2 (d {\ot} x )\, \hat{\SR}_\0(\la)$ for \hbox{$x=a,b$}.
It is easy to see that the unique (up to a scalar factor)
solution to these equations is given by the same formula~\rf{Rsol}.
But for $\theta=0$ we have $\delta(D_q)=bc \ot D_q$, which has
nontrivial operator dependence only in its first tensor component.
This makes $\hat{\SR}_\0(\la)$ rather useless for constructing
integrals of motion since it produces only those that have no
interaction between different
sites of the lattice (cf.~\rf{Hloc2} for $\theta=0$).

Thus, we see a kind of dual pictures for matrices $\g(\la)$
and $\hat{\g}(\la)$: it is the fundamental transfer--matrix
for the former and the auxiliary transfer--matrix for the latter
that generate a set of mutually commuting
elements of~$\CB^{\SRN}$.

\subsection{q--oscillator algebra ${\CA}_q$} 
\label{QOA}

Interrelations between deformed oscillator algebras and
quantum Lie algebras are well known (see, e.g.~\cite{CP,KS}).
Relation of the former to quantum groups is also known,
see e.g.~\cite{S1,DK}, but has been employed in the context
of integrable models less extensively. Here we will show that
a reduction of $\Gqe$ yields a q--oscillator algebra.
This will allow us to adapt the results of the previous sections,
in particular, the constructions of fundamental R--operators,
to the case of the q--oscillator algebra. Recall that, as above,
we deal with the case~\hbox{$q=e^{i\gamma}$}, $\gamma\in(0,\pi)$.
\pagebreak[4]

\begin{defn}\label{QA}
The q--oscillator algebra $\Aq$ is a unital associative algebra
with generators $e$, $f$, $k$, $k^{-1}$ and defining
relations $k \, k^{-1} = k^{-1} k = {\sf 1}$, and
\beq{defA}
 e \, k = q \, k\, e  \,, \qquad
 f \, k = q^{-1}\, k \, f \,, \qquad
 [e,f]=({q{-}q^{-1}}) \, k^{2} \,,
\end{equation}
and equipped with an anti--involution~* defined on generators by
\beq{A*}
  e^* = e \,, \quad f^* = f \,, \quad
  k^*=k \,, \quad (k^{-1})^* = k^{-1} \,.
\end{equation}
\end{defn}
\begin{lem}\label{cent3}
For a generic $q$, the center of $\Aq$ is generated by the
Casimir element
\beq{Cas}
 C_{q} \equiv  e \, f - q \, k^{2} \,.
\end{equation}
\end{lem}
The lemma is standard~\cite{CP,KS}.
Now we need the following simple but useful statements
which are straightforward to verify:
\begin{lem}\label{HomQ}
Let $\FU$ be $\GqeRp$ or $\GqeRpp$.
The linear homomorphism ${\CQ} : \FU \to \Aq$ defined
on generators as follows
\beq{SA}
 \CQ(a) = e \,, \quad \CQ(c) = k \,, \quad
 \CQ(b) = k\,, \quad \CQ(\theta)= k^{-1} \,, \quad \CQ(d) = f
\end{equation}
is an algebra homomorphism.
\end{lem}
\begin{lem}\label{ggAq}
The defining relations \rf{defA} are equivalent to the
following relation
\beq{RggA}
  R_{\1\2} \ \CQ(g)_{\!\1\3} \, \CQ(g)_{\!\2\3} =
  \CQ(g)_{\!\2\3} \, \CQ(g)_{\!\1\3} \ R_{\1\2} \,,
\end{equation}
where $\CQ(g)= \left(\begin{smallmatrix} e & k \\
 k & f \end{smallmatrix}\right)$,
and the auxiliary \hbox{R--matrix} is given by either of the
matrices in~\rf{Rconst}.
\end{lem}
For $\CQ$--images of the Casimir elements we have
$\CQ(D_q) = C_q$ and $\CQ(\eta'_q) = \CQ(\eta''_q) = {\sf 1}$.
The latter equalities mean that we identified $\theta$ as
the inverse to {\em both} $b$ and~$c$.

Let us introduce the following $\CQ$--images of
$g(\la)$ and $\hat{g}(\la)$:
\beq{Lqosc1}
 L^{\klein \CA}(\la) =
 \left( \!\! \begin{array}{cc}
  e & \la \,  k \\
  \la \,  k^{-1} + \la^{-1} k & f
 \end{array} \! \right) \,, \qquad
 \hat{L}^{\klein \CA}(\la) =
 \left(  \!\! \begin{array}{cc}
  \la \,  k^{-1} + \la^{-1} k & \la^{-1} f \\
  \la \, e & \la \,  k
 \end{array} \! \right) \,.
\end{equation}

\begin{propn}\label{RLA}
Each of the following matrix relations
\begin{align}
\label{RgglA}
  R_{\1\2}(\la) \ L^{\klein \CA}_{\1\3}(\la\mu) \,
  L^{\klein \CA}_{\2\3}(\mu) &=
  L^{\klein \CA}_{\2\3}(\mu) \,
 L^{\klein \CA}_{\1\3}(\la\mu) \ R_{\1\2}(\la) \,, \\
\label{RhgglA}
 \hat{R}_{\1\2}(\la) \ \hat{L}^{\klein \CA}_{\1\3}(\la\mu) \,
 \hat{L}^{\klein \CA}_{\2\3}(\mu) &=
 \hat{L}^{\klein \CA}_{\2\3}(\mu) \,
 \hat{L}^{\klein \CA}_{\1\3}(\la\mu) \ \hat{R}_{\1\2}(\la) \,,
\end{align}
where the auxiliary \hbox{R--matrices} are given by \rf{Rqosc}
and \rf{Rqosc2}, respectively, holds if and only if the elements
$e$, $f$, $k$ satisfy relations~\rf{defA} and $k^{-1}$ satisfies
the following relations:
\beq{k'}
 e \, k^{-1} = q^{-1}  k^{-1}  e \,, \qquad
 f \, k^{-1} = q \, k^{-1} \, f \,, \qquad
 [k,k^{-1}] = 0 \,.
\end{equation}
\end{propn}

\proof First, applying Lemma~\ref{HomQ} to
Eqs.~\rf{Rggl}--\rf{Rhggl}, we conclude that relations
\rf{RgglA}--\rf{RhgglA} do hold. Next, it is easy to see that all
the steps in the proof of Proposition~\ref{RGG} remain valid.
Therefore, each of relations \rf{RgglA}--\rf{RhgglA} is
equivalent to \rf{k'} together with~\rf{RggA}.
The latter matrix relation is equivalent to \rf{defA}
by Lemma~\ref{ggAq}.
\QED

Notice that the comultiplication $\Delta$ has no consistent
reduction to $\Aq$ since
$\bigl(\CQ\,{\ot}\,\CQ\bigr)\Delta(b)\neq
 \bigl(\CQ\,{\ot}\,\CQ\bigr)\Delta(c)$.
Nevertheless, it is useful to observe the following.
\begin{propn}\label{delGA}
The linear homomorphism
$\Delta_{\!\klein\CA} \colon \GqR \to \Aq^{\ot 2}$
defined on generators as follows:
$\Delta_{\!\klein\CA}(x)=\bigl(\CQ\,{\ot}\,\CQ\bigr) \Delta(x)$
for $x=a,b,c,d$, is an {\em algebra homomorphism} and a
{\em \hbox{*--homomorphism}} w.r.t. the anti--involution~\rf{A*}.
\end{propn}
\proof
The assertion follows by combining Lemma~\ref{HomQ} with the
properties of the standard comultiplication $\Delta$ for~$\Gq$.
\QED

For the non--standard ``comultiplication'' we have
the following reduction of $\delta$ to~$\Aq$.
\begin{propn}\label{delGA2}
The linear homomorphism
$\delta_{\!\klein\CA} \colon \Aq \to \Aq^{\ot 2}$
defined on generators as follows
\beq{delAq}
{}  \delta_{\!\klein\CA}(e) = e \ot k^{-1} + k \ot e \,,\qquad
 \delta_{\!\klein\CA}(f) = k \ot f   \,, \qquad
{}  \delta_{\!\klein\CA}(k^{\pm 1}) = k^{\pm 1} \ot  k^{\pm 1}
\end{equation}
is a {\em coassociative algebra homomorphism} and a
{\em \hbox{*--homomorphism}} w.r.t. the anti--involution~\rf{A*}.
\end{propn}

\proof
Notice that
\hbox{$\delta_{\!\klein\CA} \circ \CQ =(\CQ{\ot}\CQ) \circ \delta$}
is a linear homomorphism from $\FU$ to $\Aq^{\ot 2}$,
where $\FU$ is $\GqeRp$ or $\GqeRpp$.
Therefore, applying $\CQ{\ot}\CQ$ to \rf{delgh}, we infer that
\beq{delQgh}
 (id \ot \delta_{\!\klein\CA}) \, \CQ(g^\pm) =
 \CQ(g^\pm)_{\!\1\2} \, \CQ(g^\pm)_{\!\1\3} \,,
\end{equation}
where $\CQ(g^+)= \left(\begin{smallmatrix} k^{-1} & 0 \\
 e & k \end{smallmatrix}\right)$ and
$\CQ(g^-)= \left(\begin{smallmatrix} k & f \\
 0 & 0 \end{smallmatrix}\right)$.
Further we can proceed exactly as in the proof of
Proposition~\ref{Bialg3}.
\QED

\begin{defn}\label{PIa}
Let $\CB$ be an algebra of linear operators acting on a Hilbert
space~$\CK$. An irreducible representation
$\pi_{\klein \CA} \colon \Aq \to \CB$ is called {\em positive} if
the following operators are {\em self--adjoint} and
{\em strictly positive} on~$\CK$: \\
i) $\pi_{\klein \CA}(x)$ for $x=e,f,k,C_q$;\\
ii) $q^{\frac{1}{2}}\pi_{\klein \CA}(e)
 \bigl(\pi_{\klein \CA}(k)\bigr)^{-1}$\, and\,
$q^{\frac{1}{2}} \bigl(\pi_{\klein \CA}(k)\bigr)^{-1}
 \pi_{\klein \CA}(f)$.
\end{defn}

\begin{propn}\label{RA}
Let $\CB$ and $\mathcal K$ be as in Definition~\ref{PIa} and
let $\pi_{\klein \CA}$ be a positive representation of~$\Aq$.
Let $L^{\klein \CA}(\la), \hat{L}^{\klein \CA}(\la) %
 \,{\in}\, {\mathrm{Mat}}(2){\ot}\CB$ be as in~\rf{Lqosc1}.
Then the operators
$\SR^{\klein \CA}(\la), \hat{\SR}^{\klein \CA}(\la) \in \CB^{\ot 2}$
acting on ${\mathcal K}{\ot}{\mathcal K}$ and defined by the formulae
\begin{align}
\label{Rqo1}
  \SR^{\klein \CA}(\la) &= (k \ot k)^{-\frac{\al}{2} \log \la} \,
  \bigl( (e \ot k + k \ot f)(k \ot e + f \ot k) \bigr)^{\al \log \la}
   \, (k \ot k)^{-\frac{\al}{2} \log \la}  \,, \\
\label{Rqo2}
\hat{\SR}^{\klein \CA}(\la) &=
 \bigl( e k \ot k^{-1} f + k^2 \ot e f  - q k^2 \ot k^2
   \bigr)^{\al \log \la} \,,
\end{align}
where~$ \al \equiv \fr{1}{\log q}$, satisfy the equations
\begin{align}
\label{RLLA1}
 \SR^{\klein \CA}_{\2\3}(\la) \, L^{\klein \CA}_{\1\2}(\la \mu) \,
 L^{\klein \CA}_{\1\3} (\mu) &=  L^{\klein \CA}_{\1\2} (\mu) \,
 L^{\klein \CA}_{\1\3}(\la \mu) \, \SR^{\klein \CA}_{\2\3}(\la) \,,\\
\label{RLLA2}
 \hat{\SR}^{\klein \CA}_{\2\3}(\la) \,
 \hat{L}^{\klein \CA}_{\1\2}(\la \mu) \,
 \hat{L}^{\klein \CA}_{\1\3} (\mu) &=
 \hat{L}^{\klein \CA}_{\1\2} (\mu) \,
 \hat{L}^{\klein \CA}_{\1\3}(\la \mu) \,
 \hat{\SR}^{\klein \CA}_{\2\3}(\la) \,.
\end{align}
If the tensor product $\pi_{\klein \CA} \ot \pi_{\klein \CA}$
is multiplicity free, then \rf{Rqo2} is the unique solution
of \rf{RLLA2} up to multiplication by a scalar factor.
\end{propn}

\proof
First, $\pi_\0 \equiv \pi_{\klein \CA} \ot \CQ \colon \FU\to \CB$
is clearly a positive representation for $\FU=\GqeRp$ (as
well as for $\FU=\GqeRpp$). Next, it is obvious that
$\SR^{\klein \CA}(\la)$ solving \rf{RLLA1} is a solution of
Eqs.~\rf{Rcc}--\rf{Rac1}, or, equivalently,
$\tilde{\SR}^{\klein \CA}(\la) =
(k \ot k )^{\frac{\al}{2} \log \la} \,
 \SR^{\klein \CA}(\la) \, (k \ot k)^{\frac{\al}{2} \log \la}$
is a solution of Eqs.~\rf{Rcc2}--\rf{Rac2}, where each term
$x\,{\ot}\,y$ is understood as $\pi_\0(x)\,{\ot}\,\pi_\0(y)$
and $\Delta_\pi$ is replaced with
$\Delta_{\pi_\0} \equiv (\pi_{\klein \CA}\,{\ot}\,\pi_{\klein \CA})
 \circ \Delta_{\klein \CA}$.
Notice that the definition of $\Delta_{\klein \CA}$ given
in Proposition~\ref{delGA} is sufficient because
$\Delta_\pi(\theta)$ does not enter Eqs.~\rf{Rcc2}--\rf{Rac2}.
Now, it is easy to see that the $\pi_\0$ counterparts of
Eqs.~\rf{adc}-\rf{bcc} hold. This, along with Proposition~\ref{delGA},
implies that the $\pi_\0$ counterpart of formula \rf{Rl1}
holds as well. Whence we obtain formula \rf{Rqo1} as the
$\pi_\0$ counterpart of formula~\rf{Ral}. Finally, it is
easy to see that the remaining verification of Eq.~\rf{rx}
in Theorem~\ref{RG} is valid for the $\pi_\0$ counterpart
of~$\SX(\la)$.

Analogous consideration for the $\pi_\0$ counterparts
of Eqs.~\rf{SRd}--\rf{X2}, where $\delta$ (which actually
stands for $\delta_\pi$) is replaced with
$\delta_{\pi_\0} \equiv (\pi_\0 \,{\ot}\,\pi_\0) \circ \delta$,
is straightforward because, by Proposition~\ref{delGA2},
$\delta_{\pi_\0}$ has the same algebra homomorphism properties
as~$\delta_\pi$. For the same reason, the $\pi_\0$ analogue of the
uniqueness part of Theorem~\ref{RGh} is valid if we invoke
Lemma~\ref{cent3} instead of Lemma~\ref{cent2}.
\QED

Fundamental R--operators \rf{Rqo1} and \rf{Rqo2} are regular
and have the properties given in~\rf{unit}.
The corresponding local Hamiltonian densities constructed via
\rf{Hloc} are $\CQ$--images of those in \rf{Hloc1} and \rf{Hloc2},
namely
\begin{align}
\label{Hqo1}
  H^{\klein \CA}_{n,n \p1} &=
 \log \bigl( (e_{n \p1} k_{n} + k_{n \p1} f_{n})
 (k_{n \p1} e_{n} + f_{n \p1} k_{n}) \bigr) - \log(k_n k_{n \p1})
 \,,\\
\label{Hqo2}
 \hat{H}^{\klein \CA}_{n,n \p1} &=
 \log \bigl( e_{n \p1} \, k_{n \p1} k_n^{-1} f_n \,  +
  k^2_{n \p1} \, (C_q)_n   \bigr) \,.
\end{align}
As before, the arguments of the logarithms here are positive
self--adjoint operators.

\subsection{q--DST model}
\label{QDST}
The discrete self--trapping model, which describes a chain of $\SRN$
coupled anharmonic oscillators, is know to be integrable~\cite{E1,KSS}.
The corresponding L--matrix satisfies a counterpart of Eq.~\rf{RLL0}
with additive spectral parameter and the rational auxiliary R--matrix,
which is obtained {}from \rf{Rqosc} in the limit \hbox{$q\to 1$}.
It was suggested in \cite{PS,KP} that the following
\hbox{L--matrix}
\beq{Lqdst}
  L_n^{\rm\klein qDST}(\la) =
 \left( \begin{array}{cc}
  \la \,  k_n^{-1} + \la^{-1} k_n & f_n \\
  e_n & \la \,  k_n \end{array} \right) \,,
\end{equation}
where each triple ($e_n$, $f_n$, $k_n$) satisfies relations~\rf{defA}
and operators assigned to different sites commute, can be
regarded as an \hbox{L--matrix} associated with a q--deformed discrete
self--trapping (q--DST) model. The expansion of the corresponding
auxiliary transfer--matrix about the point $\la=0$ yields
\begin{eqnarray}
\label{Hqdst2}
 & T(\la) = \la^{-\SRN} \, Q  +
 \la^{2-\SRN} \, Q \cdot H^{\rm\klein qDST} + \ldots \,, &  \\
\label{Hqdst}
 &  Q = \prod\limits_{n=1}^{\SRN} k_n \,, \qquad
 H^{\rm\klein qDST} = \sum \limits_{n=1}^{\SRN}
    k^{-2}_{n} + k_n^{-1} e_n k^{-1}_{n \p1} f_{n \p1} \,. &
\end{eqnarray}
Here $Q=e^{\gamma\,h}$ with $h$ being the number of particles
operator and $H^{\rm\klein qDST}$ is a nearest--neighbour Hamiltonian
for the q--DST model.

Let us remark that \rf{Lqdst} is related to
$\hat{L}^{\klein \CA}(\la)$ in \rf{Lqosc1} via a twist
in either the auxiliary or in the quantum space:
\beq{Lh2}
 L^{\klein qDST}(\la) =
 \lambda^{\frac{1}{2}\sigma_\3} \,
 \hat{L}^{\klein \CA}(\la) \, \lambda^{-\frac{1}{2}\sigma_\3} =
 k^{\al \log \la} \, \hat{L}^{\klein \CA}(\la) \,
 k^{-\al \log \la} \,, \qquad \al=\fr{1}{\log q} \,.
\end{equation}
The first equality here implies that $L^{\klein qDST}(\la)$
satisfies Eq.~\rf{RLL0} with the same auxiliary
R--matrix~\rf{Rqosc2}. The second equality implies that the
fundamental R--operator for $L^{\klein qDST}(\la)$ is related to
that for $\hat{L}^{\klein \CA}(\la)$ as follows
(notice that \rf{Rqo2} commutes with $k\,{\ot}\,k$)
\beq{Rqdst}
 \SR^{\rm\klein qDST}(\la) =
 ({\sf 1} \ot k)^{\al \log \la} \, \hat{\SR}^{\klein \CA}(\la)\,
 (k \ot {\sf 1})^{- \al \log \la} \,.
\end{equation}
Using \rf{Rqo2} and applying formula \rf{Hloc} to \rf{Rqdst}, we find
a nearest--neighbour Hamiltonian different from \rf{Hqdst} which
corresponds to \hbox{L--matrix}~\rf{Lqdst}:
\beq{Hqdst3}
 \tilde{H}^{\rm\klein qDST} = \fr{1}{\gamma} \sum_n
  \Bigl( \log\bigl( C_q \, k^2_{n \p1} +
  e_{n \p1} k_{n \p1} k_n^{-1} f_n \bigr)
  + \log(k_n k_{n \p1}^{-1}) \Bigr) \,.
\end{equation}
Notice that the term $\log(k_n k^{-1}_{n \p1})$ does not contribute
to the total Hamiltonian in the case of a periodic chain.

\begin{rem}
Substituting \rf{Ralhdil} into \rf{Rqdst}, we obtain
(omitting a scalar factor):
\beq{Rqdst2}
 \SR^{\rm\klein qDST}(\la) =
\bigl( k \ot k  \bigr)^{\al \log \la} \
 \frac{ S_\omega\bigl(\la^{-1} \sr \bigr)}{S_\omega(\la\, \sr)} \,,
 \qquad \sr =  \bigl(C_q\bigr)^{-1} \,  k^{-1} e \ot k^{-1} f \,.
\end{equation}
An analogous formula was proposed in \cite{KP} in the case
of~$|q|<1$ in terms of the compact quantum dilogarithm~$S(x)$.
\end{rem}

\subsection{Weyl algebra} 

For the factor algebras of $\GqeRp$ and $\GqeRpp$ over the ideals
generated by the relations $c=0$ and $b=0$, respectively,
the only nontrivial defining relations are of the Weyl type.
These factor algebras are isomorphic to the following algebra.

\label{WAL}
\begin{defn}\label{WA}
The Weyl algebra $\Wq$ is a unital associative algebra with
generators $u$, $\tu$, $v$, $v^{-1}$ and defining relations
$v \, v^{-1} = v^{-1} v = {\sf 1}$ and
\beq{defW}
 u \, \tu = \tu \, u \,, \qquad u \, v = q \, v \, u  \,,\qquad
 \tu \, v = q^{-1} \, v \, \tu
\end{equation}
and equipped with an anti--involution * defined on generators by
\beq{W*}
 u^*= u \,,\qquad \tu^*= \tu\,, \qquad v^* = v \,,
 \qquad (v^{-1})^* = v^* \,.
\end{equation}
\end{defn}

The following statements are straightforward to verify.

\begin{lem}\label{cent4}
For a generic $q$, the center of $\Wq$ is generated by the
element~$Z_q = u \, \tu$.
\end{lem}
\begin{lem}\label{HomQw}
The linear homomorphisms $\CQ' : \GqeRp \to \Wq$ and
$\CQ'' : \GqeRpp \to \Wq$
defined
on generators as follows
\begin{align}
\label{QW1}
{}& \CQ'(a) = u \,, \quad \CQ'(b) = v \,, \quad \CQ'(c) = 0\,,
 \quad \CQ'(\theta)= v^{-1} \,, \quad \CQ'(d) = \tu \,, \\
\label{QW2}
{}& \CQ''(a) = u \,, \quad \CQ''(b) = 0 \,, \quad \CQ''(c) = v\,,
 \quad \CQ''(\theta)= v^{-1} \,, \quad \CQ''(d) = \tu
\end{align}
are algebra homomorphisms.
\end{lem}

Now we will introduce contractions of the maps $\Delta$ and $\delta$
suitable for~$\Wq$.

\begin{defn}\label{PIw}
Let $\CB$ be an algebra of linear operators acting on a Hilbert
space~$\CK$. An irreducible representation
$\pi_{\klein \CW} \colon \Wq \to \CB$ is called {\em positive} if
the following operators are {\em self--adjoint} and
{\em strictly positive} on~$\CK$: \\
i) $\pi_{\klein \CW}(x)$ for $x=u,\tu,v$;\\
ii) $q^{\frac{1}{2}}\pi_{\klein \CW}(u)
 \bigl(\pi_{\klein \CW}(v)\bigr)^{-1}$\, and\,
$q^{\frac{1}{2}} \bigl(\pi_{\klein \CW}(v)\bigr)^{-1}
 \pi_{\klein \CW}(\tu)$.
\end{defn}

\begin{propn}\label{Bialg4}
Let $\CB$ and $\CK$ be as in Definition~\ref{PIw}
and let $\pi_{\klein \CW}$ be a positive representation of~$\Wq$.
The linear homomorphism $\Delta_{\klein \CW} : \Wq \to \CB^{\ot 2}$
defined on generators as follows:
\beq{delw1}
\begin{aligned}
 \Delta_{\klein \CW}(u) &=
    \pi_{\klein \CW}(u) \ot \pi_{\klein \CW}(u) \,, \qquad
 \Delta_{\klein \CW}(\tu) =
    \pi_{\klein \CW}(\tu) \ot \pi_{\klein \CW}(\tu) \,, \\
 \Delta_{\klein \CW}(v) &=
    \pi_{\klein \CW}(u) \ot \pi_{\klein \CW}(v) +
    \pi_{\klein \CW}(v) \ot \pi_{\klein \CW}(\tu) \,,\qquad
 \Delta_{\klein \CW}(v^{-1}) = \bigl( \Delta_{\klein \CW}(v) \bigr)^{-1}
\end{aligned}
\end{equation}
is an {\em algebra homomorphism} and a
{\em \hbox{*--homomorphism}} w.r.t. the anti--involution~\rf{W*}.
\end{propn}
\proof
Notice that
$\Delta_{\klein \CW} \circ \CQ' =
\bigl((\pi_{\klein \CW} \circ \CQ') \ot
(\pi_{\klein \CW} \circ \CQ') \bigr) \circ \Delta$ is,
by Lemma~\ref{HomQw}, a linear homomorphism {}from $\GqeRp$ to~$\CB$.
Therefore, for $x=u,\tu,v$, the claimed properties of
$\Delta_{\klein \CW}$ are inherited from those of~$\Delta$.
For $\Delta_{\klein \CW}(v^{-1})$, a consideration analogous to that
in the proof of Proposition~\ref{Bialg2} applies since,
by Remark~\ref{Rem2}, $\Delta_{\klein \CW}(v)$ is a positive
self--adjoint and hence invertible operator.
\QED

\begin{rem}\label{DWop}
We used in this proof that $\Delta_{\klein \CW}$ is related
to $\Delta$ via~$\CQ'$. The {\em opposite} comultiplication
$\Delta'_{\klein \CW}$ (obtained by exchanging the tensor factors in
$\Delta_{\klein \CW}(x)$) is similarly related to $\Delta$ via~$\CQ''$,
namely, $\Delta'_{\klein \CW} \circ \CQ'' =
\bigl((\pi_{\klein \CW} \circ \CQ'') \ot
(\pi_{\klein \CW} \circ \CQ'') \bigr) \circ \Delta$.
\end{rem}
\begin{rem}
Using the relation between $\Delta$ and $\Delta_{\klein \CW}$,
we can write an explicit expression for $\Delta_{\klein \CW}(v^{-1})$.
Namely, applying $\CQ'$ to \rf{delthb}, we obtain
\beq{delthv}
 \Delta_{\klein \CW}(v^{-1}) =
 \bigl( \pi_{\klein \CW}\,{\ot}\,\pi_{\klein \CW} \bigr) \,\Bigl(
 S_\omega(\sw) \, \bigl(u^{-1} \,{\ot}\, v^{-1} \bigr) \,
 \bigl(S_\omega(\sw)\bigr)^{-1} \Bigr)\,,
\end{equation}
where $\sw = v u^{-1} \,{\ot}\, v^{-1}\tu $
and $\pi_{\klein \CW}(u^{-1})=\pi_{\klein \CW}(\tu)
 \bigl(\pi_{\klein \CW}(Z_q)\bigr)^{-1}$.
\end{rem}

\begin{propn}\label{Bialg5}
The linear homomorphism
$\delta_{\klein \CW} \colon \Wq \to \Wq^{\ot 2}$
defined on generators as follows:
\beq{delw2}
 \delta_{\klein \CW}(u) = u \ot v^{-1} \,, \qquad
   \delta_{\klein \CW}(\tu) = v \ot \tu \,, \qquad
 \delta_{\klein \CW}(v^{\pm 1}) =  v^{\pm 1} \ot v^{\pm 1}
 \end{equation}
is a {\em coassociative algebra homomorphism} and a
{\em \hbox{*--homomorphism}} w.r.t. the anti--involution~\rf{W*}.
\end{propn}
\proof
Straightforward. However, it is instructive to notice that
\hbox{$\delta_{\!\klein\CW} \circ \CQ'' =(\CQ''{\ot}\CQ'') \circ \delta$}
is a linear homomorphism from $\GqeRpp$ to $\Wq^{\ot 2}$.
Therefore, applying $\CQ''{\ot}\CQ''$ to \rf{delgh}, we infer that
\beq{delQgw}
 (id \ot \delta_{\!\klein\CW}) \, \CQ''(g^\pm) =
 \CQ''(g^\pm)_{\!\1\2} \, \CQ''(g^\pm)_{\!\1\3} \,,
\end{equation}
where $\CQ''(g^+)= \left(\begin{smallmatrix} v^{-1} & 0 \\
 u & 0 \end{smallmatrix}\right)$ and
$\CQ''(g^-)= \left(\begin{smallmatrix} v & \tu \\
 0 & 0 \end{smallmatrix}\right)$.
\QED

It is easy to check that any monomial in $\Wq^{\ot 2}$ which
commutes with $\delta_{\klein \CW}(x)$, $x=u,\tu,v$ is a
power of~$\delta_{\klein \CW}(Z_q)$.
But the centralizer of $\Delta_{\klein \CW}\bigl(\CW_q\bigr)$
contains not only functions of $\Delta_{\klein \CW}(Z_q)$.
\begin{lem}\label{Whw}
Denote $\sz= u v \ot u v^{-1}$. Then for all $x\in\Wq$ we have
\beq{whw}
 [(\pi_{\klein \CW} \,{\ot}\, \pi_{\klein \CW})(\sz),
    \Delta_{\klein \CW}(x)] = 0  \,.
\end{equation}
\end{lem}
\proof
It suffices to verify \rf{whw} for the generators $u,\tu,v$,
which is straightforward.
\QED


\subsection{Fundamental R--operators for $g'(\la)$ and
$\check{g}'(\la)$,  Volterra  and lattice free field models}

Let us introduce the following $\CQ'$--images of
$g(\la)$ and $\check{g}(\la)=\sigma_\1 g(\la)$:
\beq{Lvolt}
 g'(\la) =
 \left( \!\! \begin{array}{cc}
  u & \la \,  v \\
  \la \, v^{-1} & \tu
 \end{array} \! \right) \,, \qquad
 \check{g}'(\la) =
 \left(  \!\! \begin{array}{cc}
  \la \, v^{-1} & \tu \\
   u & \la \, v
 \end{array} \! \right) \,.
\end{equation}
Matrix $g'(\la)$ is the L--matrix for the Volterra model~\cite{V2}
and is also related to the lattice sine--Gordon model~\cite{V2,F1,F2}.
We will see below that $\check{g}'(\la)$ is the L--matrix for
the Volterra model for a dual dynamical variable (we use $\check{g}'(\la)$
rather than $\CQ'(\hat{g}(\la))$ to make the duality most transparent;
the corresponding fundamental R--operators differ only by a twist).
In the compact case, a fundamental R--operator for $g'(\la)$ was found
in~\cite{V2}. Here we will give an alternative derivation, which exhibits
transparently the underlying comultiplication structure.
For brevity of notations, we will write $x\,{\ot}\,y$
instead of $\pi_{\klein \CW}(x)\,{\ot}\,\pi_{\klein \CW}(y)$.

\begin{thm}\label{Rg'}
Let $\CB$ and $\mathcal K$ be as in Definition~\ref{PIw} and
let $\pi_{\klein \CW}$ be a positive representation of~$\Wq$.
Let $g'(\la), \check{g}'(\la) \,{\in}\, {\mathrm{Mat}}(2){\ot}\CB$
be as in~\rf{Lvolt}. Then the operators
$\SR'(\la), \check{\SR}'(\la) \in \CB^{\ot 2}$
acting on ${\mathcal K}{\ot}{\mathcal K}$ and defined by the formulae
\begin{align}
\label{Rv1}
  \SR'(\la) &= r(\sz,\la) \ \tilde{\sz}^{\frac{\al}{4} \log \la} \,
  \bigl( u \ot v + v \ot \tu \bigr)^{\al \log \la}
   \, \tilde{\sz}^{\frac{\al}{4} \log \la}  \,, \\
\label{Rv2}
\check{\SR}'(\la) &=
  r(\tilde{\sz},\la) \
 \sz^{\frac{\al}{4} \log \la} \,
  \bigl( \tu \ot v + v^{-1} \ot \tu \bigr)^{\al \log \la}
   \, \sz^{\frac{\al}{4} \log \la}  \,,
\end{align}
where $\sz= u v \,{\ot}\, u v^{-1}$,
$\tilde{\sz}=\tu v^{-1} \,{\ot}\, u v^{-1}$,
$\al \equiv \fr{1}{\log q}$, satisfy the equations
\begin{align}
\label{RLLV1}
 \SR'_{\2\3}(\la) \, g'_{\1\2}(\la \mu) \,
 g'_{\1\3} (\mu) &=  g'_{\1\2} (\mu) \,
 g'_{\1\3}(\la \mu) \, \SR'_{\2\3}(\la) \,,\\
\label{RLLV2}
 \check{\SR}'_{\2\3}(\la) \,
 \check{g}'_{\1\2}(\la \mu) \, \check{g}'_{\1\3} (\mu) &=
 \check{g}'_{\1\2} (\mu) \, \check{g}'_{\1\3}(\la \mu) \,
 \check{\SR}'_{\2\3}(\la) \,,
\end{align}
for any choice of the function $r(t,\la)$.
\end{thm}

\proof
Eq.~\rf{RLLV1} can be regarded as a $\CQ'{\ot}\CQ'$--image
of \rf{RLL1}. It is easy to see, that Eqs.~\rf{Rcc}--\rf{Rac1}
turn into the following relations:
\begin{align}
\label{Ruv0}
{} [\SR'(\la), v \ot v^{-1}] = [\SR'(\la), v^{-1} \ot v] &=
 [\SR'(\la), u \ot u] = [\SR'(\la), \tu \ot \tu] = 0 \,,\\
 \label{Ruv1}
{} \SR'(\la) \, ( u \ot v + \la \, v \ot \tu)  &=
 (\la \, u \ot v + v \ot \tu) \, \SR'(\la)\,, \\
 \label{Ruv2}
 {} \SR'(\la) \, ( \la\, v^{-1} \ot u + \tu \ot v^{-1})  &=
 (v^{-1} \ot u + \la \, \tu \ot v^{-1}) \, \SR'(\la)\,.
\end{align}
To exhibit maximally the structure of these equations related to
the comultiplication $\Delta_{\klein \CW}$, we introduce
$\tilde{\SR}'(\la) = \tilde{\sz}^{-\frac{\al}{4} \log \la} \,
 \SR'(\la) \, \tilde{\sz}^{-\frac{\al}{4} \log \la}$.
Then equations \rf{Ruv0}--\rf{Ruv2} acquire the following form:
\begin{align}
\label{Ruv3}
{} \tilde{\SR}'(\la) \, (v \ot v^{-1}) =
  \la\, (v \ot v^{-1}) \, \tilde{\SR}'(\la) \,, &\qquad
  \tilde{\SR}'(\la) \, \Delta_{\klein \CW}(v)
   = \Delta_{\klein \CW}(v) \, \tilde{\SR}'(\la) \,,\\
 \label{Ruv4}
{} \tilde{\SR}'(\la) \, \Delta_{\klein \CW}(u)
   = \la^{-1} \Delta_{\klein \CW}(u) \, \tilde{\SR}'(\la) \,, & \qquad
   \tilde{\SR}'(\la) \, \Delta_{\klein \CW}(\tu)
   = \la\, \Delta_{\klein \CW}(\tu) \, \tilde{\SR}'(\la) \,.
\end{align}

It is now easy to see that \rf{Ruv3}--\rf{Ruv4} are solved by
\beq{Rtv1}
\tilde{\SR}'(\la) = r(\sz,\la) \
    \bigl(\Delta_{\klein \CW}(v) \bigr)^{\al \log\la} \,,
\end{equation}
where $r(t,\la)$ can be an arbitrary function thanks to Lemma~\ref{Whw}
and the fact that \hbox{$[\sz, v \,{\ot}\,v^{-1}]=0$}.
Thus, we established that \rf{Rv1} satisfies \rf{Ruv0}--\rf{Ruv1}.
Taking into account that
$[\sz,v^{-1} \,{\ot}\, u]=[\sz,\tu \,{\ot}\, v^{-1}]=0$,
it remains to prove that $\SR'_\0(\sw,\la)=
\tilde{\sz}^{\frac{\al}{4} \log \la}
\bigl(\Delta_{\klein \CW}(v) \bigr)^{\al \log\la}
\tilde{\sz}^{\frac{\al}{4} \log \la}$
satisfies~\rf{Ruv2}. For this purpose we apply Lemma~\ref{LUL}
and rewrite it as follows:
\beq{Rv3}
 \SR'_\0(\sw,\la) = \bigl( Z_q \bigr)^{\al\log \la} \,
 \frac{ S_\omega\bigl(\la^{-1} \sw\bigr)}%
  {S_\omega\bigl(\la\,\sw\bigr)} \, \sw^{-\frac{\al}{2}\log\la} \,,
   \qquad \sw = v u^{-1} \ot v^{-1} \tu \,,
\end{equation}
where we used that $\sw \, \tilde{\sz} = q^4 \tilde{\sz} \, \sw$.
Multiplying \rf{Ruv2} with $v \,{\ot}\, \tu$ from the right,
we obtain the following functional equation on function
$\SR'_\0(w,\la)$:
\beq{Rw}
 \SR'_\0( w, \la)\,(\la + q^{-1} w) =
 (1 + \la\,  q^{-1} w ) \, \SR'_\0(q^{-2} w, \la) \,,
\end{equation}
which is easy to verify using Eq.~\rf{Sdil}.

To prove that \rf{Rv2} satisfies \rf{RLLV2}, we  observe that
$g'(\la)$ and $\check{g}'(\la)$ are related in a way which fits
the hypotheses of Lemma~\ref{AUT} (namely,
$s=\left(\begin{smallmatrix} 0 & 1 \\ 1 & 0 \end{smallmatrix}\right)$
and the automorphism $\iota$ is defined by
$\iota(u)=\tu$, $\iota(\tu)=u$, $\iota(v^{\pm 1})=v^{\mp 1}$).
Therefore, according to Eq.~\rf{autR},
the fundamental R--operator for $\check{g}'(\la)$
is $\check{\SR}'(\la)=(\iota \,{\ot}\, id) \SR'(\la)$.
Noticing that $\iota(\tilde{\sz})=\sz$ and
$\iota(\sz)=\tilde{\sz}$, we obtain formula~\rf{Rv2}.
\QED

\begin{rem}
In \cite{V2,F2,FV2}, another solution to Eq.~\rf{Rw} was
given, namely
\beq{Rv5}
 \tilde{\SR}'_\0( \sw, \la) =
 \frac{S_\omega(\sw) \, S_\omega(\sw^{-1}) }%
 {S_\omega(\la\, \sw) \, S_\omega(\la\, \sw^{-1})} \,.
\end{equation}
Eq.~\rf{ssw} in Appendix~\ref{ApA} shows that \rf{Rv3} and
\rf{Rv5} coincide up to a factor independent of~$\sw$.
\end{rem}

Fundamental R--operators \rf{Rv1} and \rf{Rv2} are regular
the sense of Eq.~\rf{reg} if $r(t,1)=1$. Furthermore, they
have the properties given in~\rf{unit} provided that
$\bar{r}(t,\la)=r(t,\la^{-1})=1/r(t,\la)$ for $t,\la\,{>}\,0$
(notice that $\sz^*=\sz$,\, $\tilde{\sz}^*=\tilde{\sz}$, and
$[\sz,\tilde{\sz}]=0$).
The corresponding local Hamiltonian densities constructed via
\rf{Hloc} are given by ($r'(t)$ stands for the derivative of
$r(t,\la)$ w.r.t.~$\la$ at $\la \,{=}\,0$)
\begin{align}
\label{Hv1}
 \gamma\, H'_{n,n \p1} &=
 \log \bigl( v_n u_{n \p1} + \tu_n v_{n \p1}\bigr)
  - \fr{1}{2} \log(v_n v_{n \p1}) + r'(\sz_{n \p1, n})
 + \fr{1}{2} \log(u_n \tu_{n \p1} ) \,,\\
\label{Hv2}
 \gamma\,  \check{H}'_{n,n \p1} &=
 \log \bigl( v_n \tu_{n \p1} +  \tu_n v^{-1}_{n \p1} \bigr) +
 \fr{1}{2} \log(u_n u_{n \p1} ) + r'(\tilde{\sz}_{n \p1, n}) +
 \fr{1}{2} \log(v^{-1}_n v_{n \p1})  \,.
\end{align}
The arguments of the logarithms here are positive
self--adjoint operators. Notice that the last terms in
\rf{Hv1} and \rf{Hv2} add only a constant to the total
Hamiltonian in the case of a periodic chain.

Consider the following  positive representations of $\Wq$
on the Hilbert space $\CK = L^2(\BR)$
\begin{align}
\label{EFK2a}
{} \pi_+(u) &= e^{ p} \,,\qquad
 \pi_+(\tu) = e^{ -p} \,,\qquad
 \pi_+(v) = e^{- 2\phi} \,, \\
\label{EFK2b}
{} \pi_-(u) &= e^{ -2\phi} \,,\qquad
 \pi_-(\tu) = e^{2\phi} \,,\qquad
 \pi_-(v) = e^{- p} \,,
\end{align}
were  $p$ and $\phi$ are self--adjoint operators which satisfy
$[p,\phi]=\fr{\gamma}{2i}$, $\gamma\in(0,\pi)$. For these
representations, the classical limit of \rf{Hv1}--\rf{Hv2}
acquires the following form (up to additive constants):
\begin{align}
\label{Hv1b}
 \gamma\, \pi_+(H'_{n,n \p1}) &=
 \log \cosh s_+    + r'(e^{2s_-}) \,,\\
\label{Hv2b}
 \gamma\,  \pi_-(\check{H}'_{n,n \p1} )
{} &= \log \cosh s_-   + r'(e^{2s_+}) \,.
\end{align}
where $s_\pm \equiv \fr{1}{2} p_n \,{+}\, \fr{1}{2} p_{n \p1}
  \,{\pm}\,  \phi_{n \p1} \,{\mp}\, \phi_n$.
It was shown in \cite{V2} that $s_\pm$ are related
(in a nonultralocal way, via a discretized Miura transformation)
to the dual dynamical variables of the Volterra model, and that
\rf{Hv1b} for $r(t,\la)=0$ coincides with the Hamiltonian of
the Volterra model for~$s_+$. The obvious symmetry between
\rf{Hv1b} and \rf{Hv2b} makes it clear that the matrix
$\check{g}'(\la)$ can as well be taken as an L--matrix for
the Volterra model and that the corresponding Hamiltonian
\rf{Hv2b} for $r(t,\la)=0$ is the Hamiltonian of the Volterra
model for the dual dynamical variable~$s_-$.

Let us demonstrate that $g'(\la)$  can be regarded also as
an L--matrix for a lattice regularization of the free field.
For this goal we have to choose such $r(t,\la)$ in \rf{Rv1}
that $r'(e^{2t})=\log \cosh t \,{+}\, {\rm const}$ in
the classical limit. For instance, we can take
(cf.~\rf{Rv3} and notice that $[\sz,\sw]=0$ and
$\sz_{n,n \p 1} \sw_{n \p 1,n} = Z_q$)
\beq{Rff1}
 \SR'(\la) =
  \frac{ S_\omega\bigl(\la^{-1} Z_q \, \sz^{-1}\bigr)}%
  {S_\omega\bigl(\la\,Z_q\,\sz^{-1}\bigr)} \,
  \bigl(\sz\, \sw^{-1}\bigr)^{\frac{\al}{2}\log\la} \,
\frac{ S_\omega\bigl(\la^{-1} \sw\bigr)}%
  {S_\omega\bigl(\la\,\sw\bigr)} \,.
\end{equation}
Then \rf{Hv1b} acquires the following form:
\beq{Hv1c}
 \gamma\, \pi_+(H'_{n,n \p1}) =
 \log \cosh s_+    +  \log \cosh s_- \,.
\end{equation}
In the continuum limit \rf{cont}, we have
$s_\pm = \kappa \bigl(p(x) \pm \partial_x \phi(x)\bigr) + o(\kappa)$
($\kappa$ stands for the lattice spacing) and \rf{Hv1c} turns
into $H'_{n,n \p1}={\rm const} +
\frac{\kappa^2}{\gamma}\bigl(p^2 + (\partial_x \phi)^2\bigr)
+ o(\kappa^2)$, i.e., it becomes the Hamiltonian density
of the free field. Furthermore, assining a copy of
$L^{\rm f}(\la) = \pi_+\bigl(g'(\kappa\la)\bigr)$
to each site of the lattice, we get the following continuum
limit of this L--matrix: $ L_n^{\rm f}(\la)=
 \bigl(\begin{smallmatrix} 1 & 0\\ 0 & 1
    \end{smallmatrix} \bigr)
 + \kappa \, \bigl(U_+(\la)+U_-(\la)\bigr) + O(\kappa^2)$,
where
\beq{UccL2}
  U_+(\la) =   \left( \begin{array}{cc}
  \fr{1}{2} \, \partial_+ \phi  & \la \, e^{ - 2\phi } \\
  \la \, e^{2\phi}  & - \fr{1}{2} \, \partial_+ \phi
 \end{array} \right) \,, \quad
 U_-(\la) =   \left( \begin{array}{cc}
  \fr{1}{2} \, \partial_- \phi  &  0  \\
  0 &  - \fr{1}{2} \, \partial_- \phi
 \end{array} \right) \,.
\end{equation}
These matrices satisfy the zero curvature equation,
$\partial_-U_+(\la) + \partial_+U_-(\la) = 2[U_+(\la),U_-(\la)]$,
provided that $\phi$ satisfies the equation of motion of the free
field: $\square\, \phi = 0$.

\subsection{Fundamental R--operator for $g''(\la)$,
lattice free field model}

Let us introduce the following $\CQ''$--image of $g(\la)$:
\beq{Lff}
 g''(\la) =
 \left( \!\! \begin{array}{cc}
  u &  0 \\
  \la \, v^{-1} + \la^{-1}  v& \tu
 \end{array} \! \right) \,.
\end{equation}

\begin{thm}\label{Rg''}
Let $\CB$ and $\mathcal K$ be as in Definition~\ref{PIw} and
let $\pi_{\klein \CW}$ be a positive representation of~$\Wq$.
Let $g''(\la) \,{\in}\, {\mathrm{Mat}}(2){\ot}\CB$
be as in~\rf{Lff}. Then the operator $\SR''(\la)  \in \CB^{\ot 2}$
acting on ${\mathcal K}{\ot}{\mathcal K}$ and defined by the formula
\beq{Rff2}
  \SR''(\la) = \hat{\sz}^{\frac{\al}{2} \log \la} \,
 \Bigl( \bigl( v \ot u + \tu \ot v \bigr)\,
 \bigl( u \ot v + v \ot \tu \bigr) \Bigr)^{2 \al \log \la} \,
 \hat{\sz}^{\frac{\al}{2} \log \la} \, \,,
\end{equation}
where $\hat{\sz}=\tu v^{-2} \,{\ot}\, u v^{-2}$ and
$\al \equiv \fr{1}{\log q}$, satisfies the equation
\beq{RLLff}
 \SR''_{\2\3}(\la) \, g''_{\1\2}(\la \mu) \,
 g''_{\1\3} (\mu) =  g''_{\1\2} (\mu) \,
 g''_{\1\3}(\la \mu) \, \SR''_{\2\3}(\la) \,.
\end{equation}
\end{thm}

\proof
As always, for brevity of notations, we write $x\,{\ot}\,y$
instead of $\pi_{\klein \CW}(x)\,{\ot}\,\pi_{\klein \CW}(y)$.
Eq.~\rf{RLLff} can be regarded as a $\CQ''{\ot}\CQ''$--image
of~\rf{RLL1}. It is easy to see, that Eqs.~\rf{Rcc}--\rf{Rac1}
turn into the following relations:
\begin{align}
\label{Rab0}
{} \SR''(\la) \, ( v \ot u + \la \, \tu \ot v)  &=
 (\la \, v \ot u + \tu \ot v) \, \SR''(\la)\,, \\
 \label{Raa0}
{} [ \SR''(\la) , u \ot u ] =0 \,, \quad & \qquad
 [ \SR''(\la) , \tu \ot \tu ] =0 \,, \\
{} \SR''(\la) \, (\la\, v^{-1} \ot u + \tu \ot v^{-1} )  &=
 (v^{-1} \ot u + \la\, \tu \ot v^{-1}) \, \SR''(\la) \,.
 \label{Rac0}
\end{align}
It is easy to recognize in \rf{Rab0}--\rf{Raa0} a structure
related to the opposite comultiplication~$\Delta'_{\klein \CW}$
in accordance with Remark~\ref{DWop}.
To make this structure more transparent, we introduce
\beq{tsr''}
\tilde{\SR}''(\la) = (v \ot v )^{\frac{\al}{2} \log \la} \,
 \SR''(\la) \, (v \ot v)^{\frac{\al}{2} \log \la} \,.
\end{equation}
Then equations \rf{Rab0}--\rf{Rac0} acquire the following form:
\begin{align}
 \label{Rvp}
 \tilde{\SR}''(\la) \, \Delta'_{\klein \CW}(v) &=
 \Delta'_{\klein \CW}(v) \, \tilde{\SR}''(\la) \,, \\
 \label{Rad0}
 \tilde{\SR}''(\la) \, \Delta'_{\klein \CW}(u) =
 \la^{-2} \Delta'_{\klein \CW}(u) \, \tilde{\SR}''(\la) \,,
  \quad & \qquad
 \tilde{\SR}''(\la) \, \Delta'_{\klein \CW}(\tu) =
 \la^2 \, \Delta'_{\klein \CW}(\tu) \, \tilde{\SR}''(\la) \,, \\
\label{Rac00}
 \tilde{\SR}''(\la) \, (\la\, v^{-1} \ot u + \la^{-1} \tu \ot v^{-1} )
 &= (\la^{-1} v^{-1} \ot u + \la\, \tu \ot v^{-1}) \,
    \tilde{\SR}''(\la) \,.
\end{align}
According to Lemma~\ref{Whw}, a solution to Eqs.~\rf{Rvp}--\rf{Rad0}
may contain as a factor an arbitrary function of
$\check{\sz}= u v^{-1} \ot u v$. Actually, it is more convenient to
introduce
$\sw \equiv Z_q \,\check{\sz}^{-1}=v u^{-1} \,{\ot}\, v^{-1} \tu$.
Then \rf{Rvp}--\rf{Rad0} are solved by
\beq{Rw2}
\tilde{\SR}''(\la) =
 \bigl(\Delta'_{\klein \CW}(v)\bigr)^{2 \al \log \la}
 \, \check{\SR}''(\sw,\la) \,,
\end{equation}
where $\check{\SR}''(t,\la)$ is yet undetermined function. Noticing that
$[\Delta'_{\klein \CW}(v),v^{-1} \,{\ot}\, u]=
 [\Delta'_{\klein \CW}(v),\tu \,{\ot}\, v^{-1}]=0$, we infer that
$\check{\SR}''(\sw,\la)$ must solve~\rf{Rac00}.
Multiplying \rf{Rac00} with $u \,{\ot}\, v$ {}from the right,
we obtain the following functional equation on function
$\check{\SR}''(w, \la)$:
\beq{Rw3}
 \check{\SR}''( w, \la)\,(\la \, q w^{-1} + \la^{-1}) =
 (\la^{-1}  q w^{-1} + \la) \, \check{\SR}''(q^{-2} w, \la) \,.
\end{equation}
Comparing this equation with \rf{Rtv1}--\rf{Rw}, we
conclude that
\beq{Rv4}
 \check{\SR}''(\sw,\la) = \bigl( Z_q \bigr)^{2\al\log\la}
 \frac{ S_\omega\bigl(\la^{-2} \sw\bigr)}%
  {S_\omega\bigl(\la^2 \sw\bigr)} \, \sw^{-\al\log\la}
 = \tilde{\sz}^{\frac{\al}{2} \log \la}
\bigl(\Delta_{\klein \CW}(v) \bigr)^{2\al \log\la}
\tilde{\sz}^{\frac{\al}{2} \log \la} \,,
\end{equation}
where $\tilde{\sz}=\tu v^{-1} \,{\ot}\, u v^{-1}$.
Notice that $\Delta'_{\klein \CW}(v)$ commutes with
$\Delta_{\klein \CW}(v)$ and~$\tilde{\sz}$.
Combining \rf{tsr''}, \rf{Rw2}, and \rf{Rv4}, we
obtain formula~\rf{Rff2}.
\QED

Fundamental R--operator \rf{Rff2} is regular and has the
properties given in~\rf{unit}. The corresponding local
Hamiltonian density constructed via \rf{Hloc} is given by
\beq{Hff1}
 {\gamma} \, H''_{n,n \p1} =
 2 \log \Bigl( \bigl( v_{n \p1} u_n + \tu_{n \p1} v_n \bigr)\,
 \bigl( u_{n \p1} v_n + v_{n \p1} \tu_n \bigr) \Bigr)  +
 \log(\tu_{n \p1} v^{-2}_{n \p1} u_n v^{-2}_n) \,.
\end{equation}
Definition~\ref{PIw} along with Remark~\ref{Rem2} ensure that
the argument of the first logarithm here is a product of commuting
positive self--adjoint operators. In the classical limit \rf{Hff1}
can be written as follows:
\beq{Hff1c}
 {\gamma} \, H^{\prime\prime, \rm cl}_{n,n \p1} =
 2 \log \Bigl( u_n u_{n \p1} + \tu_n \tu_{n \p1} +
 Z_q \bigl( v_n^{-1} v_{n \p1} + v_n v^{-1}_{n \p1} \bigr) \Bigr)
 + \log(\tu_{n \p1}  u_n ) \,.
\end{equation}

Consider the following one--parameter family $\pi_{\kappa}$ of
positive representations of~$\Wq$:
\beq{EFK3}
 \pi_{\kappa}(u)  = \fr{1}{\kappa} \, e^{\frac{\beta}{4}\, \Pi} \,\qquad
 \pi_{\kappa}(\tu)  = \fr{1}{\kappa} \, e^{-\frac{\beta}{4}\, \Pi} \,,
 \qquad \pi_{\kappa}(v) = e^{- \frac{\beta}{2} \, \Phi}\,,
\end{equation}
where $\kappa$, $\beta$, $\Pi$, and  $\Phi$ are as in \rf{EFK1}.
Let us introduce the following \hbox{L--matrix}:
$ L^{\klein\rm F}(\la) = \kappa\, \pi_\kappa \bigl(g''(\la) \bigr)$
and assign a copy of this matrix to each site of the lattice,
\beq{ff}
 L_n^{\klein \rm F}(\la)   =
 \left( \begin{array}{cc}
  e^{\frac{\beta}{4}\, \Pi_n}   &   0 \\
  \kappa\, \bigl(\la \, e^{\frac{\beta}{2} \Phi_n} +
  \la^{-1}  e^{-\frac{\beta}{2} \Phi_n} \bigr) &
   e^{- \frac{\beta}{4} \, \Pi_n}
 \end{array} \right) \,,
\end{equation}
where \hbox{$[\Pi_n,\Phi_m]=-i \delta_{nm}$}.
L--matrix \rf{ff} can be obtained from the Liouville L--matrix
if in Eq.~\rf{LL} we shift the zero mode of the field:
$\Phi_n \to \Phi_n + \xi$, rescale the spectral parameter:
$\la\to \la\, e^{-\xi \frac{\beta}{2}}$, and take the limit
$\xi\to +\infty$.

One may expect that in such a limit the Liouville model
turns into the free field. Indeed, for the representation
\rf{EFK3}, Eq.~\rf{Hff1c} acquires the following form:
\beq{Hff1d}
 {\gamma} \, H^{\prime\prime, \rm {\klein F}, cl}_{n,n \p1} =
 2 \log \bigl( 2 \cosh \fr{\beta}{4}(\Pi_n \,{+}\, \Pi_{n \p1}) +
 2 \cosh \fr{\beta}{2}(\Phi_{n \p1} \,{-}\, \Phi_n)  \bigr)
 + {\rm const}\,,
\end{equation}
where we omitted the last term in \rf{Hff1c} since it does not
contribute to the total Hamiltonian in the case of a periodic chain.
In the continuum limit \rf{cont}, we recover from \rf{Hff1d}
the Hamiltonian density of the free field:
$H^{\prime\prime, \rm {\klein F}, cl}_{n,n \p1}={\rm const} +
 \kappa^2 \bigl(\Pi^2 + (\partial_x \Phi)^2\bigr) + o(\kappa^2)$.
Furthermore, in the continuum limit we have
$L_n^{\klein \rm F}(\la) =
 \bigl(\begin{smallmatrix} 1 & 0\\ 0 & 1
    \end{smallmatrix} \bigr)
 + \kappa \, \bigl(U_+(\la)+U_-(\la)\bigr) + O(\kappa^2)$,
where
\beq{UccL3}
  U_+(\la) =   \left( \begin{array}{cc}
  \fr{\beta}{8} \, \partial_+ \Phi  & 0 \\
   \la \, e^{\frac{\beta}{2} \Phi_n}  &
 - \fr{\beta}{8} \, \partial_+ \Phi
 \end{array} \right) \,, \quad
 U_-(\la) =   \left( \begin{array}{cc}
  \fr{\beta}{8} \, \partial_- \Phi  &  0  \\
  \la^{-1}  e^{-\frac{\beta}{2} \Phi_n} &
 - \fr{\beta}{8} \, \partial_- \Phi
 \end{array} \right) \,.
\end{equation}
These matrices satisfy the zero curvature equation,
$\partial_-U_+(\la) + \partial_+U_-(\la) = 2[U_+(\la),U_-(\la)]$,
provided that $\Phi$ satisfies the equation of motion of the free
field: $\square\, \Phi = 0$.

\rem
Let us remark that the two fundamental R--operators that we have
found for the lattice free field are quite similar.
Namely, it is straightfoward to check that
\beq{fF}
 \SR''(\la) = (u^{-1} {\ot}\, u)^{\frac{\al}{2} \log \la} \
 \SR'(\la^2) \ (u^{-1} {\ot}\, u)^{\frac{\al}{2} \log \la} \,,
\end{equation}
where $\SR''(\la)$ is given by \rf{Rff2} and
$\SR'(\la)$ is given by \rf{Rff1}. Notice that, for a periodic chain,
the factors $(u^{-1} {\ot} u)^{\frac{\al}{2} \log \la}$ do not
contribute to the total Hamiltonian.

\subsection{Fundamental R--operator for $\hat{g}''(\la)$,
relativistic Toda model}

Let us introduce the following $\CQ''$--image of $\hat{g}(\la)$:
\beq{Lto}
 \hat{g}''(\la) =
 \left( \!\! \begin{array}{cc}
  \la \, v^{-1} + \la^{-1}  v& \la^{-1}  \tu \\
  \la \, u &  0
 \end{array} \! \right) \,.
\end{equation}
This matrix is related via a twist (cf.~\rf{Lh2}) to the L--matrix
of the relativistic Toda model~\cite{KT,PS}:
\beq{LToda}
  L^{\rm\klein rT}(\la) = (\pi_- \,{\ot}\, \pi_-)
 \bigl(\fr{1}{i} v^{\al \log \la} \,
    \hat{g}''(i\la) \, v^{-\al \log \la} \bigr) =
 \left( \begin{array}{cc}
  \la \, e^{p} - \la^{-1} e^{-p} & -e^{2\phi} \\
  e^{-2\phi} & 0 \end{array} \right) \,,
\end{equation}
where $\pi_-$ is the positive representation \rf{EFK2b} of~$\Wq$.
A suitable limit of \rf{LToda} for \hbox{$q\to 1$} yields the
L--matrix of the ordinary Toda chain model, which satisfies a counterpart
of Eq.~\rf{RLL0} with additive spectral parameter and a rational
auxiliary R--matrix.

Integrals of motion both for the ordinary and relativistic Toda models
are constructed by means of expanding the auxiliary transfer--matrix
$T(\la)$ (cf.~Section~\ref{QDST}). Results of the present article
explain why the corresponding fundamental R--operators cannot be
employed for this purpose.

\begin{thm}
Let $\CB$ and $\mathcal K$ be as in Definition~\ref{PIw} and
let $\pi_{\klein \CW}$ be a positive representation of~$\Wq$.
Let $\hat{g}''(\la) \,{\in}\, {\mathrm{Mat}}(2){\ot}\CB$
be as in~\rf{Lto}. Then the operator $\hat{\SR}''(\la)\in \CB^{\ot 2}$
acting on ${\mathcal K}{\ot}{\mathcal K}$ and defined by the formula
\beq{Rto1}
 \hat{\SR}''(\la) =
 \bigl(\delta_{\klein \CW} (Z_q) \bigr)^{\al \log\la} =
 (u \, v \ot v^{-1} \tu)^{\al \log\la} \,,
\end{equation}
where $\al \equiv \fr{1}{\log q}$, satisfies the equations
\beq{RLLT1}
 \hat{\SR}''_{\2\3}(\la) \, \hat{g}''_{\1\2}(\la \mu) \,
 \hat{g}''_{\1\3} (\mu) =  \hat{g}''_{\1\2} (\mu) \,
 \hat{g}''_{\1\3}(\la \mu) \, \hat{\SR}''_{\2\3}(\la) \,.
\end{equation}
\end{thm}

\proof
Reexamining the proof of Theorem~\ref{RGh} in the case of $b=0$,
we see that an analogue of Eq.~\rf{SRd} holds in the form
$[\hat{\SR}''(\la),\delta_{\klein \CW} (x)]$, $x=u,\tu,v$.
Unlike $\Delta_{\klein \CW}$, the centralizer of
$\delta_{\klein \CW}(x)$, $x\in\Wq$ is generated only
by~$\delta_{\klein \CW}(Z_q)$.
Therefore, $\hat{\SR}''(\la)$ must be a function of
$\delta_{\klein \CW} (Z_q)$. It is easy to see that that
the $b=0$ counterparts of Eqs.~\rf{Rabd}--\rf{Re2} determine
this function uniquely (up to a scalar factor) and lead to
formula~\rf{Rto1}.
\QED

Fundamental R--operator \rf{Rto1} is regular and has the properties
given in~\rf{unit}. However, the corresponding local Hamiltonian
density constructed via \rf{Hloc},
\beq{Hto}
 \gamma\, H''_{n,n \p1} =
 \log \bigl( v^{-1}_n \tu_n u_{n \p1} v_{n \p1}\bigr) \,,
\end{equation}
leads to a trivial total Hamiltonian in the case of a periodic chain.

\section{Conclusion}

We have developed the Baxterization approach to the quantum
group $\Gq$ and emphasized the role of the standard
and non--standard comultiplications for constructing
the corresponding fundamental R--operators.
Our results imply that the quantum symmetry algebra for a number
of integrable lattice models is the quantum group $\Gq$ or its
reductions for which the comultiplication structure is a reduction
of those for~$\Gq$. This is especially remarkable in the
case of the lattice Liouville model because the quantum
group $\Gq$ itself emerged for the first time exactly in the
study of relations for the monodromy matrix of the lattice
Liouville model~\cite{FT1}.
For the Volterra model, we have shown that the two
dual L--matrices lead to the same Hamiltonian but for
the dual dynamical variables. We have also emphasized the
role of the ambiguity in the solution for the corresponding
fundamental R--operators: fixing it in a trivial way yields
the Hamiltonian of the Volterra model, whereas fixing it
in a self--dual way yields the Hamiltonian of a lattice
regularization of the free field. For the latter model
we have also found another L--matrix which can
be regarded as a limit of that for the lattice
Liouville model. It is interesting that, although
the free field in continuum is a very simple model,
the fundamental R--operators related to its lattice
regularization have quite a nontrivial structure.

\par\vspace*{3mm}\noindent
{\small
{\bf Acknowledgement.}
 This work was supported in part by
 Alexander von Humboldt Foundation, by the Russian Science
 Support Foundation, and by the Russian Foundation for
 Fundamental Research (grants 07--02--92166 and 08--01--00638).
 The author is grateful to V.~Schomerus and J.~Teschner
 for useful discussions and for warm hospitality
 at DESY (Hamburg), where a part of this work was done.
 The author thanks K.~Schm\"udgen for useful remarks.
}

\appendix

\section{Appendix}

\subsection{Quantum dilogarithm}\label{ApA}

Consider the functional equation
\beq{Sdil}
  S(q^{-1} x) = (1+x) \, S(q\, x)   \,.
\end{equation}
Its solution is given by the product
$S(x)=\prod_{n=1}^{\infty}(1+ x\,q^{2n-1})$,
which is convergent for \hbox{$|q|<1$}. This function appears in various
related forms in lattice integrable models~\cite{T1,V2,FV1} and was
coined ``quantum dilogarithm'' in~\cite{FK1}. It was observed in
\cite{F2,F3} that, for $q=e^{i \pi \omega^2}$, \hbox{$\omega\in (0,1)$},
a well--defined solution to \rf{Sdil} is given by
\beq{qdil}
  S_\omega(x)=\prod_{n=1}^{\infty} \frac{(1+ x\,q^{2n-1})}%
 {(1+ x^{\omega^{-2}}\,\hat{q}^{2n-1})} =
 \exp \biggl\{
 \int\limits_{\Omega} \frac{dt}{4\, t} \,
 \frac{ e^{ \frac{t}{i\pi \omega} \, \log x}}%
 {\sinh \omega t \, \sinh{\frac{t}{\omega}} } \biggr\} \,,
\end{equation}
where $\hat{q} \equiv e^{-i \pi \omega^{-2}}$ and $\Omega=\BR+i0$.
Among the important properties of $S_\omega(x)$ are
\begin{align}
 \text{ self--duality:} & \qquad
 S_\omega(x^{\omega}) = S_{\omega^{-1}}(x^{\omega^{-1}}) \,, \\
\label{Sunit}
 \text{ unitarity:} & \qquad
 \overline{S_\omega(x)} \ S_\omega(x) =1 \quad
 \text{for $x\in {\mathbb R}_+$} \,.
\end{align}
This function is closely related to the Barnes double gamma function
\cite{B1} and plays an important role in studies of non--compact
quantum groups \cite{F4,PT,W1,S2,BT1,W2,V3} and related integrable
models~\cite{KLS,FK2,T2,BT2}.

The following lemma proves to be useful for converting
powers of coproducts in formulae for fundamental R--operators
into expressions involving quantum dilogarithms.
\begin{lem}\label{LUL}
Let $u$ and $v$ be a pair of positive self--adjoint
operators satisfying, in the sense of Remark~1, the Weyl relation:
$u\, v = q^2 \, v \, u$, where $q=e^{i\pi \omega^2}$, $\omega \in (0,1)$.
Suppose that $\sw \equiv q\, u^{-1} v$ is positive self--adjoint.
Then the following identity holds:
\beq{uvdil}
 (u+v)^{t }  = u^{\frac{t }{2}} \,
   \frac{S_\omega(q^{-t} \sw)}{S_\omega(q^t\, \sw)} \, u^{\frac{t}{2}}
 = v^{\frac{t}{2} } \,
 \frac{S_\omega(q^{-t} {\sw}^{-1})}{S_\omega(q^{t}\, {\sw}^{-1})} \,
  v^{\frac{t }{2} } \,.
\end{equation}
\end{lem}
\proof
Using relations $u\, f(\sw) = f(q^2 \sw) \,u$\ and
$v\, f(\sw^{-1}) = f(q^{-2} \sw^{-1}) \,v$, it is easy to verify
the following identities
\beq{uvS1}
 u+v = S_\omega(\sw) \,u\, \bigl(S_\omega(\sw)\bigr)^{-1}
 =  \bigl(S_\omega(\sw^{-1})\bigr)^{-1} \,v\, S_\omega(\sw^{-1})  \,.
\end{equation}
These are relations of unitary equivalence, thanks to the
property~\rf{Sunit}. Whence we infer that
\beq{uvS2}
 (u+v)^{t } =
 S_\omega(\sw) \,u^{\frac{t }{2}}\,u^{\frac{t }{2}}\,
 \bigl(S_\omega(\sw)\bigr)^{-1} = u^{\frac{t }{2}}\,
 \frac{S_\omega(q^{-t } \sw)}{S_\omega(q^{t } \sw)} \,
  u^{\frac{t }{2}}\,,
\end{equation}
holds in the sense of Remark~1.
The second equality in \rf{uvdil} can be derived analogously.
\QED

\begin{rem}\label{RemW}
Equality of the two expressions involving quantum dilogarithms in
\rf{uvdil} allows us to obtain the following functional identities:
\beq{ssw}
 w^t = \frac{S_\omega(q^{-t } w) \, S_\omega(q^{t} w^{-1})}%
 {S_\omega(q^{t } w) \, S_\omega(q^{-t} w^{-1}) } \ =\,
 q^{t^2} \ \frac{S_\omega(q^{-2t } w) \, S_\omega(q^{2t} w^{-1})}%
 {S_\omega(w) \, S_\omega(w^{-1}) } \,.
\end{equation}
\end{rem}

For the proof of Theorem~2 we will need the following lemma (for
the sake of brevity we will write $x$ instead of $\pi_{\CA}(x)$).
\begin{lem}\label{LG}
Let $e$, $f$, and $k$ generate a positive representation
$\pi_{\cal A}$ of the q--oscillator algebra~${\cal A}_q$
(cf. Definition~\ref{PIa}). Then the following relation holds
\beq{rbd3}
 \frac{ S_\omega(\la^{-1} f )}{S_\omega\bigl(\la\, f \bigr)} \,
 \bigl(\la\, k^2 + e\bigr) = \bigl(\la^{-1} k^2 + e \bigr) \,
 \frac{ S_\omega(\la^{-1} f )}%
    {S_\omega\bigl(\la\, f \bigr)}  \,.
\end{equation}
\end{lem}
For a positive representation, we can write $e=u_e +v_e$,
where $u_e= f^{-1} C_q$ and $v_e=q^{-1}f^{-1}k^2$ are positive
self--adjoint operators satisfying relation
$u_e v_e =q^2v_e u_e$ (hence, by Remark~\ref{Rem2}, $e$~is
positive self--adjoint; a rigorous operator--theoretic
consideration of the formula $e=f^{-1}(C_q + q^{-1}k^2)$
is given in~\cite{S1}).
Therefore, if $G(t)$ is a sufficiently nice function (i.e.
$G(f)$ has a suitable domain, cf. the discussion in~\cite{S1}),
then we have $k^2 G(f) = G(q^2 f) k^2$ and
$e G(f) = G(f) C_q f^{-1} + q^{-1} G(q^2 f) f^{-1} k^2$.
Taking these relations into account, we infer that the
operator equation
\beq{rbd3'}
 G(f,\la) \, \bigl(\la\, k^2 + e\bigr)
 = \bigl(\la^{-1} k^2 + e \bigr) \, G(f,\la)
\end{equation}
is equivalent to the following functional one:
\beq{rbd3''}
 G(f,\la) \, \bigl(\la + q^{-1} f^{-1} \bigr)
 = \bigl(\la^{-1} + q^{-1} f^{-1} \bigr) \, G(q^2 f,\la) \,.
\end{equation}
Using \rf{Sdil}, it is straightforward to check that
$G(f,\la) = \dfrac{ S_\omega(\la^{-1} f )}%
    {S_\omega\bigl(\la\, f \bigr)}$ solves~\rf{rbd3''}.
\QED

\subsection{Proof of Lemma~\ref{RBD}}\label{ApB}
Formula \rf{Ralhdil} can be rewritten as follows:
\beq{Ralh'}
  \hat{\SR}(\la) = \bigl( D_q \bigr)^{\al \log \la} \,
 \bigl( b \ot {\sf 1} \bigr)^{\al \log \la} \, \check{\SR}(\sr;\la) \,
 \bigl( c \ot {\sf 1} \bigr)^{\al \log \la} \,,
\end{equation}
where
\beq{Ralh''}
  \check{\SR}(\sr;\la) =
 \frac{ S_\omega(\la^{-1} \sr )}{S_\omega(\la\, \sr )}
\,,\qquad \sr = \bigl(D_q\bigr)^{-1} \,
 \bigl( q^{-\frac{1}{2}} b^{-1} a \bigr) \ot
 \bigl( q^{\frac{1}{2}} \theta d \bigr) \,.
\end{equation}
Substituting \rf{Ralh'} in \rf{Re1}, it is easy to check
that Lemma~\ref{RBD} is equivalent to the assertion that
$\check{\SR}(r;\la)$ satisfies the following relation:
\beq{rbd2}
 \check{\SR}(\sr;\la) \, (\la\, \theta \ot d + d \ot b) =
 (\la^{-1} \theta \ot d + d \ot b) \, \check{\SR}(\sr;\la) \,.
\end{equation}

Notice that $d \ot b$, $\theta \ot d$, and $\sr$  are positive
self--adjoint operators. Now, a simple computation (using, in
particular, the identity $q\,da\,{-}\,q^{-1}ad=(q\,{-}\,q^{-1})D_q$)
yields
\begin{eqnarray}
\label{rhodb}
 & (d \ot b) \, \sr - \sr \, (d \ot b) =
    (q-q^{-1}) \, \theta \ot d  & \\
\label{dbtd}
 & (d \ot b) \, (\theta \ot d) =
    q^2 \, (\theta \ot d) \, (d \ot b) \,, \qquad
 \sr \, (\theta \ot d) = q^{-2} (\theta \ot d) \, \sr \,. &
\end{eqnarray}
Comparing these relations with \rf{defA} we see that
$\hat{e} = d \ot b$, $\hat{f} = \sr$, and
$\hat{k}^2 = \theta \ot d$ generate a positive representation
of the algebra~${\cal A}_q$ ($\hat{k}$ can be defined as the
unique positive self--adjoint square root of $\theta \ot d$).
Invoking Lemma~\ref{LG}, we establish validity of Eq.~\rf{rbd2}
and hence of Lemma~\ref{RBD}.

\def\baselinestretch{1}

\end{document}